\documentclass[11pt,reqno]{amsart}
\usepackage{amssymb}
\usepackage{enumerate}
\usepackage{calc}

\theoremstyle{plain}
\newtheorem{thm}{Theorem}[section]
\newtheorem{lem}[thm]{Lemma}
\newtheorem{cor}[thm]{Corollary}

\newtheorem{prop}[thm]{Proposition}

\theoremstyle{definition}

\newtheorem{defn}[thm]{Definition}
\newtheorem{desc}[thm]{Description}

\theoremstyle{remark}
\newtheorem{rem}[thm]{Remark}

\advance \topmargin by -\headheight
\advance \topmargin by -\headsep
     
\evensidemargin \oddsidemargin
\newenvironment{entry}
   {\begin{list}{}%
          {%
            \setlength{\labelwidth}{35pt}%
            \setlength{\topsep}{3\smallskipamount}%
            \setlength{\parsep}{2\smallskipamount}%
            \setlength{\leftmargin}{\labelwidth+\labelsep}%
          }%
   }%
   {\end{list}}

\newcommand{\w}{\ensuremath{\omega}}
\newcommand{\wa}{\ensuremath{\omega^{\alpha}}}
\renewcommand{\a}{\ensuremath{\alpha}}
\renewcommand{\b}{\ensuremath{\beta}}
\newcommand{\g}{\ensuremath{\gamma}}
\renewcommand{\epsilon}{\varepsilon}
\renewcommand{\phi}{\varphi}
\newcommand{\e}{\ensuremath{\varepsilon}}

\newcommand{\ab}[2]{\mbox{$#1\!\cdot\!#2$}}
\newcommand{\ib}[1]{I_b(#1)}
\newcommand{\llo}{\ensuremath{\ell_1}}
\newcommand{\lli}{\mbox{$\ell_1$-index}}
\newcommand{\llis}{\mbox{$\ell_1$-indices}}
\newcommand{\llt}{\mbox{$\ell_1$-tree}}
\newcommand{\llkt}{\mbox{$\ell_1$-$K$-tree}}
\newcommand{\llbbt}{\mbox{$\ell_1$-block} basis tree}

\newcommand{\llkbbt}{\mbox{$\ell_1$-$K$-block} basis tree}
\newcommand{\llp}{\mbox{$\ell_1^+$}}
\newcommand{\llpi}{\mbox{$\ell_1^+$-index}}

\newcommand{\llpt}{\mbox{$\ell_1^+$-tree}}
\newcommand{\llpk}{\mbox{$\ell_1^+$-$K$}}
\newcommand{\llpks}{\mbox{$\ell_1^+$-$K$-sequence}}
\newcommand{\llpkt}{\mbox{$\ell_1^+$-$K$-tree}}
\newcommand{\llpbbt}{\mbox{$\ell_1^+$-block} basis tree}
\newcommand{\llpbbi}{\mbox{$\ell_1^+$-block} basis index}
\newcommand{\llpbbis}{\mbox{$\ell_1^+$-block} basis indices}
\newcommand{\llpkbbt}{\mbox{$\ell_1^+$-$K$-block} basis tree}

\newcommand{\llpwi}{\mbox{$\ell_1^+$-weakly null index}}

\newcommand{\llpwt}{\mbox{$\ell_1^+$-weakly null tree}}

\newcommand{\llpkwt}{\mbox{$\ell_1^+$-$K$-weakly null tree}}

\newcommand{\lii}{\mbox{$\ell_{\infty}$-index}}

\newcommand{\lip}{\mbox{$\ell_{\infty}^+$}}
\newcommand{\lipi}{\mbox{$\ell_{\infty}^+$-index}}

\newcommand{\lipt}{\mbox{$\ell_{\infty}^+$-tree}}
\newcommand{\lipk}{\mbox{$\ell_{\infty}^+$-$K$}}

\newcommand{\lipbbi}{\mbox{$\ell_{\infty}^+$-block} basis index}

\newcommand{\wstar}{\mbox{weak$^*$}}
\newcommand{\ws}{\mbox{wide-$(s)$}}
\newcommand{\wss}{\mbox{wide-$(s)$ sequence}}
\newcommand{\wst}{\mbox{wide-$(s)$ tree}}
\newcommand{\lws}{\mbox{$\lambda$-wide-$(s)$}}

\newcommand{\mws}{\mbox{$\mu$-wide-$(s)$}} 
 
\newcommand{\wc}{\mbox{wide-$(c)$}}
\newcommand{\wcs}{\mbox{wide-$(c)$ sequence}}

\newcommand{\lwc}{\mbox{$\lambda$-wide-$(c)$}}

\newcommand{\mwc}{\mbox{$\mu$-wide-$(c)$}} 
 
\newcommand{\sst}{\mbox{$s$-node}}
\newcommand{\sss}{\mbox{$s$-subsequence}}

\newcommand{\spp}{\operatorname{supp}}

\newcommand{\diam}{\operatorname{diam}}
\newcommand{\sign}{\operatorname{sign}}

\newcommand{\R}{\ensuremath{\mathbf{R}}}
\newcommand{\N}{\ensuremath{\mathbf{N}}}
\newcommand{\cb}{\ensuremath{\mathcal{B}}}

\newcommand{\cs}[1][\alpha]{\ensuremath{\mathcal{S}_{#1}}}

\newcommand{\kequiv}[1][K]{\stackrel{#1}{\sim}\mbox{uvb }}

\newcommand{\C}[1][\alpha]{\ensuremath{C(\omega^{\omega^{#1}})}}
\newcommand{\Co}[1][\alpha]{\ensuremath{C_0(\omega^{\omega^{#1}})}}
\newcommand{\B}[1][\alpha]{\ensuremath{\cb(\omega^{\omega^{#1}})}}

\newcommand{\CC}[2][\alpha]{ %
  \ensuremath{C(\omega^{\omega^{#1}\cdot{#2}})}}
\newcommand{\CCo}[2][\alpha]{ %
  \ensuremath{C_0(\omega^{\omega^{#1}\cdot{#2}})}}

\newcommand{\frc}[2][1]{\ensuremath{\frac{#1}{#2}}}
\newcommand{\wkstar}{\stackrel{w^*}{\rightarrow}}
\newcommand{\wkly}{\stackrel{w}{\rightarrow}}
\renewcommand{\P}[1][\varepsilon]{\ensuremath{\mathbf{P}(#1)}}
\newcommand{\Q}[1][\varepsilon]{\ensuremath{\mathbf{Q}(#1)}}
\newcommand{\ch}{\mathbf{1}}
\newcommand{\prq}{\preccurlyeq}
\newcommand{\chir}{\text{\raisebox{\depth}{\large $\chi$}}}
\newcommand{\zetar}{\text{\raisebox{0.5\depth}{$\zeta$}}}
\newcommand{\psif}[1]{\psi{\scriptscriptstyle #1}}
\newcommand{\seq}[2][1]{\ensuremath{(#2)_{#1}^{\infty}}}
\newcommand{\io}{\ensuremath{\iota}}
\newcommand{\sgn}{\operatorname{sign}}
\newcommand{\xim}{\ensuremath{(x_i)_1^m}}
\newcommand{\yin}{\ensuremath{(y_i)_1^n}}
\newcommand{\yjn}{\ensuremath{(y_j)_1^n}}
\newcommand{\xmyn}{\ensuremath{(x_1,\dots,x_m,y_1,\dots,y_n)}}

\begin{document}

\title[The Szlenk index and local $\ell_1$-indices]{The Szlenk index and
  local $\boldsymbol{\ell_1}$-indices}
\author{Dale Alspach}
\address{Department of Mathematics\\
  Oklahoma State University\\
  Stillwater, OK 74078-1058 \\
  U.S.A. } 
\email{alspach@math.okstate.edu} 
\author[Robert Judd]{Robert Judd}
\address{Department of Mathematics\\
  University of Missouri-Columbia\\
  Columbia, MO 65211 \\
  U.S.A. } 
\email{rjudd@math.missouri.edu} 
\author[Edward Odell]{Edward Odell}
\address{Department of Mathematics
  University of Texas at Austin\\
  Austin, TX 78712-1082 \\
  U.S.A. }
\email{odell@math.utexas.edu}
\thanks{The last named author was supported by the NSF and TARP}

\begin{abstract}
  We introduce two new local $\ell_1$-indices of the same type as the
  Bourgain \lli; the $\ell_1^+$-index and the $\ell_1^+$-weakly null
  index.  We show that the $\ell_1^+$-weakly null index of a Banach space
  $X$ is the same as the Szlenk index of $X$, provided $X$ does not
  contain $\ell_1$.  The $\ell_1^+$-weakly null index has the same form as
  the Bourgain \lli: if it is countable it must take values $\w^{\a}$ for
  some $\a<\w_1$.  The different \llis\ are closely related and so knowing
  the Szlenk index of a Banach space helps us calculate its \lli, via the
  $\ell_1^+$-weakly null index.  We show that $I(\C)=\w^{1+\a+1}$.
\end{abstract}
\maketitle

\section{Introduction}
\label{sec:intro}

If $X$ is a separable Banach space, then one can study the complexity of
the \llo\ substructure of $X$ via Bourgain's \llo\ ordinal index
$I(X)$,~\cite{b} (defined precisely below).  One has $I(X)<\omega_1$ if
and 
only if \llo\ does not embed into $X$.  It was shown in~\cite{jo} that
$I(X)=\w^{\a}$ for some $\a<\w_1$ provided $\llo\not\hookrightarrow X$.
If $X$ has a basis, then one can also define an \llo\ block basis index
$I_b(X)$,~\cite{jo}.  In this paper we introduce and study five additional
related isomorphically invariant indices: $I^+(X)$, $I^+_b(X)$, $J^+(X)$,
$J^+_b(X)$ and $I^+_w(X)$.  The latter we call the \llpwi\ and show it is
equal to the Szlenk index of $X$ provided that \llo\ does not embed into
$X$.  The \llpi\ $I^+(X)$, and \llpbbis\ are motivated by the fundamental
work of James~\cite{j2}, and of Milman and Milman~\cite{mm}, on bases and
reflexivity.  These results yield that the \llpi\ is countable if and
only if $X$ is reflexive, and is equal to \w\ if and only if $X$ is
super-reflexive.  The \llpbbi\ measures the ``shrinkingness'' of a basis.
The \lipi, $J^+(X)$, and the \lipbbi\ are the obvious dual notions to the
\llp-indices, and the \lipbbi\ measures the ``boundedly completeness'' of
a basis.

All the indices are defined in terms of certain trees on $X$.  We give the
necessary background on trees in Section~\ref{sec:trees} and define the
indices in Section~\ref{sec:l1indices}.  In that section we also obtain a
number of results concerning these indices.  In Section~\ref{sec:szlenk}
we recall the Szlenk index and discuss its relation with the \llpwi.
Section~\ref{sec:Xa-Ca-indices} is concerned with calculating the various
indices for two particular collections of Banach space: the $C(\a)$ spaces
and the generalized Schreier spaces $X_{\a}$ for $\a<\w_1$.

\section{General Trees}
\label{sec:trees}

In this section we review the basic definitions and properties of the
trees we will be using.  We then construct certain specific trees.  These
trees will be abstract sets and may be thought of as ``tree skeletons.''
The nodes don't have any meaning on their own; they merely serve as a
frame on which to hang our Banach space trees.

\begin{defn}\label{defn:trees}
  By a \emph{tree} we shall mean a non-empty, partially ordered set
  $(T,\leq)$ for which the set $\{ y\in T : y\leq x\}$ is linearly ordered
  and finite for each $x\in T$.  The elements of $T$ are called
  \emph{nodes}.  The \emph{predecessor node} of $x$ is the maximal element
  $x'$ of the set $\{ y\in T : y<x\}$, so that if $y<x$, then $y\leq x'$.
  An \emph{immediate successor} of $x\in T$ is any node $y>x$ such that
  $x\leq z\leq y$ implies that $z=x$ or $z=y$.  The \emph{initial nodes}
  of $T$ are the minimal elements of $T$ and the \emph{terminal nodes} are
  the maximal elements.  A \emph{branch} of a tree is a maximal linearly
  ordered subset of a tree.  A \emph{subtree} of a tree $T$ is a subset of
  $T$ with the induced ordering from $T$.  This is clearly again a tree.
  Further, if $T'\subset T$ is a subtree of $T$ and $x\in T$, then we
  write $x<T'$ to mean $x<y$ for every $ y\in T'$.  We will also consider
  trees related to some fixed set $X$.  \emph{A tree on a set $X$} is a
  partially ordered subset $T\subseteq\cup_{n=1}^{\infty}X^n$ such that
  for $(x_1,\dots,x_m), (y_1,\dots,y_n) \in T$,
  $(x_1,\dots,x_m)\leq(y_1,\dots,y_n)$ if and only if $m\leq n$ and
  $x_i=y_i$ for $i=1,\dots,m$.
\end{defn}

We next recall the notion of the order of a tree.  Let the \emph{derived
  tree} of a tree $T$ be $D(T)=\{ x\in T : x<y \mbox{ for some } y\in
T\}$.  It is easy to see that this is simply $T$ with all of its terminal
nodes removed.  We then associate a new tree $T^{\alpha}$ to each ordinal
$\alpha$ inductively as follows.  Let $T^0=T$, then given $T^{\alpha}$ let
$T^{\alpha+1}=D(T^{\alpha})$.  If $\alpha$ is a limit ordinal, and we have
defined $T^{\beta}$ for all $\beta<\alpha$, let
$T^{\alpha}=\cap_{\beta<\alpha}T^{\beta}$.  A tree $T$ is
\emph{well-founded} provided there exists no subset $S\subseteq T$ with
$S$ linearly ordered and infinite.  The \emph{order} of a tree $T$ is
defined as $o(T)=\inf\{\alpha : T^{\alpha}=\emptyset \}$ if there exists
$\a<\w_1$ with $T^{\alpha}=\emptyset$, and $o(T)=\w_1$ otherwise.

A tree $T$ on a topological space $X$ is said to be \emph{closed}
provided the set $T\cap X^n$ is closed in $X^n$, endowed with the
product topology, for each $n\geq1$.  We have the following result
(see~\cite{d}) concerning the order of a closed tree on a
Polish space.
\begin{prop}\label{prop:wf=>cble}
   If $T$ is a well-founded, closed tree on a Polish
  (separable, complete, metrizable) space, then $o(T)<\omega_1$.
\end{prop}

A map $f:T\rightarrow T'$ between trees $T$ and $T'$ is a \emph{tree
  isomorphism} if $f$ is one to one, onto and an order isomorphism
($\alpha<\beta$ if and only if $f(\alpha)<f(\beta)$).  We will write
$T\simeq T'$ if $T$ is tree isomorphic to $T'$ and
$f:T\stackrel{\sim}{\rightarrow}T'$ to denote a tree isomorphism which,
for brevity, we shall simply call an \emph{isomorphism}.

\begin{defn}\label{defn:min-tree}
  \cite{jo}~For an ordinal $\alpha<\omega_1$ a tree $S$ is a \emph{minimal
    tree of order $\alpha$} if for each tree $T$ of order $\alpha$ there
  exists a subtree $T'\subset T$ of order $\alpha$ which is isomorphic to
  $S$.  It is easy to see that if $T$ is a tree of order $\beta$, with
  $\alpha\leq\beta<\omega_1$, then there exists a subtree $T'\subseteq T$
  which is a minimal tree of order $\alpha$.  In~\cite{jo} certain minimal
  trees $T_{\a}$ for each ordinal $\alpha<\omega_1$ were constructed
  inductively as follows.  The smallest tree $T_1$ is just a single node.
  Given $T_{\alpha}$ one chooses $z\not\in T_{\alpha}$ and puts this as
  the initial element of the tree to give $T_{\alpha+1}$.  Thus
  $T_{\alpha+1}=T_{\alpha}\cup\{z\}$ with $z<x$ for every $x\in
  T_{\alpha+1}\setminus\{z\}$.  If $\alpha$ is a limit ordinal and
  $T_{\beta}$ has been constructed for each $\beta<\alpha$, then one
  chooses a sequence of ordinals $\alpha_n$ increasing to $\alpha$, and
  sets $T_{\alpha}$ to be the disjoint union of the trees $T_{\alpha_n}$.
\end{defn}

\begin{defn}\label{defn:restricted}
  Let $T$ be a tree on a set $X$ and let $T'$ be a subtree of $T$.  We
  define another tree on $X$, the \emph{restricted subtree $R(T')$ of $T'$
    with respect to $T$}.  Let $z=(x_i)_1^n\in T'$ and let $y$ be the
  unique initial node of $T'$ such that $y\leq z$; let $m\leq n$ be such
  that $y=(x_i)_1^m$.  If $y$ is also an initial node of $T$, then set
  $k=0$, otherwise let $k<m$ be such that $(x_i)_1^k$ is the predecessor
  node of $y$ in $T$.  Finally, setting $R(z)=(x_{k+1},\dots,x_n)$, we
  define $R(T')=\{ R(z) : z\in T'\}$.  It is easy to see that $o(T')\leq
  o(R(T'))$. 
\end{defn}

Many of the proofs of the results we obtain rely on extracting certain
subtrees from the trees we are given.  To do this we construct a type of
tree called a \emph{replacement tree}.  The idea is that given two trees
$S$ and $S'$, one can, in some sense, replace each node of $S$ with a tree
isomorphic to $S'$ to obtain a much larger tree.  We know that if a tree
is isomorphic to this larger tree, for some pair $S$ and $S'$, then it is
easy to reverse the replacement process and obtain a subtree isomorphic to
$S$.  We discuss two specific types of replacement tree here, $T(\a,\b)$
for $\a,\b<\w_1$ and $T(\a,s)$ for $\alpha<\omega_1$ where $s$ is the tree
which is just a countably infinite sequence of incomparable nodes.  

\begin{desc}\label{defn:rep-tree-ab}
\cite{jo}~The replacement trees $T(\alpha,\beta)$ satisfy the following
properties for each pair $\alpha,\beta<\omega_1$:
\renewcommand{\labelenumi}{(\alph{enumi})}
\renewcommand{\labelenumii}{(\roman{enumii})}
\begin{enumerate}
\item There exists a map $f_{\alpha,\beta}: T(\alpha,\beta)\rightarrow
  T_{\alpha}$ satisfying:
    \begin{enumerate}
    \item For each $x\in T_{\alpha}$ there exists $I=\{1\}$ or $\N$
      and trees $t(x,j)\simeq T_{\beta}\ (j\in I)$ such that
      $f_{\alpha,\beta}^{-1}(x)=\cup_{j\in I}t(x,j)$ (incomparable
      union) with $I=\{1\}$ if $\alpha$ is a successor ordinal and $x$
      is the unique initial node, or $\beta<\omega$, and $I=\N$
      otherwise.
    \item For each pair $a,b\in T(\alpha,\beta),\ a\leq b$ implies
      $f_{\alpha,\beta}(a)\leq f_{\alpha,\beta}(b)$.
    \end{enumerate}
  \item $o(T(\alpha,\beta))=\beta\!\cdot\!\alpha$.
  \item $T(\alpha,\beta)$ is a minimal tree of order
    $\beta\!\cdot\!\alpha$.
\end{enumerate}
The full details of the construction of these trees may be found
in~\cite{jo}.
\end{desc} 

\begin{desc}\label{defn:rep-tree-as}
  The trees $T(\a,s)$ are built up in a similar way to the minimal trees
  $T_{\alpha}$ except that at each stage an infinite sequence of nodes is
  added instead of a single node.  Let $s=\{z^1,z^2,\dots\}$, an infinite
  sequence of incomparable nodes, and then let $T(1,s)=s$.  To construct
  $T(\alpha+1,s)$ from $T(\a,s)$ we take the set $s$ and then after each
  element put a tree isomorphic to $T(\alpha,s)$.  For example, $T(n,s)$
  is a countably infinitely branching tree of $n$ levels.  If $\alpha$ is
  a limit ordinal, and we have constructed $T(\beta,s)$ for each
  $\beta<\alpha$, then we take a \emph{sequence of successor ordinals}
  $\alpha_n\nearrow\alpha$ and let $T(\alpha,s)$ be the disjoint union
  over $n\in\N$ of trees isomorphic to $T(\alpha_n,s)$.
\end{desc} 

Each tree $T(\alpha,s)$ has the following properties:
\begin{enumerate}
\item $o(T(\alpha,s))=\alpha$;
\item $T(\alpha,s)$ has an infinite sequence of initial nodes;
\item if $z$ is in the derived tree $D(T(\alpha,s))$ (i.e.\ $z$ is not a
  terminal node of $T(\alpha,s)$), then $z$ has an infinite sequence of
  immediate successors.
\end{enumerate}
If $S$ is either the sequence of immediate successors of some node $z\in
T(\alpha,s)$, so that $S=\{w\in T:z<w \text{ and } z\leq y\leq w\text{
  implies } y=z\text{ or } y=w \}$, or the sequence of initial nodes, then
we say that $S$ is an \emph{$s$-node of $T(\alpha,s)$}.  In order to use
the trees $T(\alpha,s)$ we must build in one more property; we need to put
an ordering on the \sst s.  Thus, to each \sst, $S$ of $T(\alpha,s)$, we
associate a bijection $\psi=\psi_S:\N\rightarrow S$ and then we may write
$S=\{ z^i : i\geq1 \}$, where $z^i=\psi(i)\in S$.

Let $T$ be a tree on a Banach space $X$.  When we say \emph{$T$ is
  isomorphic to $T(\a,s)$} we shall mean not only are they isomorphic as
trees, but we shall also require that if $(x_1,\dots,x_k)\in T$, then
$(x_1,\dots,x_{k-1})\in T$.  If $T$ is such a tree and $S=\{ z^i : i\geq1
\}$ is an \sst\ of $T$, with $z^i=(x_1,\dots,x_k,y_i)$, or $z^i=(y_i)$,
for each $i\geq1$, then we say that $\{ y_i : i\geq1 \}$ is an \emph{\sss\ 
  of $T$}.

\begin{defn}\label{defn:wnt}
  A \emph{weakly null tree} on a Banach space $X$ is a tree $T$ isomorphic
  to $T(\alpha,s)$, for some $\alpha<\omega_1$, such that every
  $s$-subsequence is weakly null.
\end{defn}

In many of the proofs that follow we take certain subtrees of trees
isomorphic to $T(\alpha,s)$ on $X$ for some $\alpha$.  Given such a tree
on $X$ we assume that the sequences of nodes down a branch, and the
sequences of nodes in the \sst s satisfy some property \Q, for $\e>0$.  We
take subtrees by extracting subsequences of nodes going down branches and
subsequences of the \sst s simultaneously so that these subsequences all
satisfy some property \P.

The basic idea is straightforward: given a sequence $(x_i)_1^{\infty}$
with property \Q\ we attempt to extract a subsequence
$(x_{n_i})_1^{\infty}$ with property \P.  For example, \Q\ might be the
property that the sequence is \textit{normalized and weakly null} (with no
dependence on $\e$ here), while \P\ could be the property that the
subsequence $(x_{n_i})_1^{\infty}$ is \textit{an \e-perturbation of a
  normalized block basis of a given basis $(e_i)_1^{\infty}$}, i.e.\ there
exists a normalized block basis $(b_i)_1^{\infty}$ of $(e_i)_1^{\infty}$
such that $\sum_i\|b_i-x_{n_i}\|<\e$.  Of course, given a normalized
weakly null sequence $(x_i)_1^{\infty}$ in a Banach space with basis
$(e_i)_1^{\infty}$ we can always extract a subsequence
$(x_{n_i})_1^{\infty}$ that is an \e-perturbation of a normalized block
basis of $(e_i)_1^{\infty}$.  The trick is to do this for all sequences in
a tree.

We use the same technique each time so to avoid repeating it in each proof
we present the framework below for arbitrary properties \P\ and \Q.  Let
$\phi:\N\rightarrow\N$ be given by
\[ \phi(i+n(n-1)/2)=i, \text{ where } 1\leq n \text{ and } 1\leq i\leq n \
, \] 
thus $(\phi(n))_1^{\infty} = (1,1,2,1,2,3,1,2,3,4,\dots)$. 

\begin{lem}[Pruning Lemma]\label{lem:pruning}
  Let $X$ be a Banach space, and for $\epsilon>0$ let \P\ and \Q\ be
  properties which a sequence $(x_i)$ (finite or infinite) in $X$ may
  possess, satisfying for every finite (or empty) sequence $(u_i)_1^l$
  with property \P, and for each $\delta>0$:
  \begin{enumerate}[\rm PL(1)]
  \item for all sequences $(x_i)_1^{\infty}$ in $X$ satisfying property
    \Q, there exists a subsequence $(x'_i)_1^{\infty}$ of
    $(x_i)_1^{\infty}$ such that $(u_1,\dots,u_l,x'_1,x'_2,\dots)$ has
    property \P[\epsilon+\delta] (we then say that
    \emph{$(x'_i)_1^{\infty}$ has property \P[\epsilon+\delta] for
      $(u_i)_1^l$});
  \item if $(y_{n,i})_{i=1}^{\infty}$ are sequences in $X$ ($n\geq1$)
    satisfying \Q\ and such that $(u_1,\dots,u_l,y_{n,1},y_{n,2},\dots)$
    has property \P\ for each $n\geq1$, then there exist sequences
    $(y_i)_1^{\infty}\subseteq\{ y_{n,i} : n,i\geq1 \}$ and $1\leq
    k_1<k_2<\cdots$ with $y_i=y_{\phi(i),k_i}$ and such that
    $(u_1,\dots,u_l,y_1,y_2,\dots)$ has property \P[\epsilon+\delta].
  \item if $(x_i)_1^{\infty}$ has \P, then $(x_i)_1^k$ has \P\ for every
    $k\geq1$. 
  \end{enumerate}
  Then for any $\epsilon, \delta>0$, for every finite sequence $(u_i)_1^l$
  with property \P, for every $\alpha<\omega_1$, and for every tree $T$ on
  $X$ isomorphic to $T(\alpha,s)$, if every \sss\ of $T$ satisfies \Q,
  then there exists a subtree $S$ of $T$ which is also isomorphic to
  $T(\alpha,s)$, and such that for all nodes $z=(x_i)_1^k\in S$ with
  immediate successors $z^j=(x_1,\dots,x_k,y_j)$, the sequence
  $(u_1,\dots,u_l,x_1,\dots,x_k,y_1,y_2,\dots)$ has property
  \P[\epsilon+\delta], and the sequence $(u_1,\dots,u_l,w_1,w_2,\dots)$
  has property \P[\epsilon+\delta], where $z^j=(w_j)$, $j\geq1$, are the
  initial nodes (where the nodes $z^j$ are ordered as an \sst\ of $S$).
\end{lem}

\begin{rem}\label{rem:pl}
  We sum up the conclusion of the Pruning Lemma by saying that \emph{$S$
    has property \P[\epsilon+\delta] for $(u_i)_1^l$}, and if $S$ has
  \P[\epsilon+\delta]\ for the empty sequence, then we just say that $S$
  has \P[\epsilon+\delta].
\end{rem}

\begin{proof}
  We use induction on $\alpha$; the case $\alpha=1$ follows directly from
  hypothesis PL(1).  Suppose the result is true for $\alpha$, and
  fix $\epsilon, \delta>0$, $(u_i)_1^l$ with property \P\ and let $T$ be a
  tree on $X$ isomorphic to $T(\alpha+1,s)$ such that every \sss\ of $T$
  satisfies \Q.  Let $(z^i)_1^{\infty}$ be the sequence of initial nodes
  of $T$ with $z^i=(w_i)$ for some $w_i\in X$.  Using PL(1) we may find a
  subsequence $(w'_i)_1^{\infty}$ of $(w_i)_1^{\infty}$ such that
  $(u_1,\dots,u_l, w'_1,w'_2,\dots)$ has property \P[\epsilon+\delta/2].
  Let $\bar z^i=(w'_i)$ and set $\bar T=\{ x\in T : x\geq\bar z^i \mbox{
    for some } i \}$.  This tree is still isomorphic to $T(\alpha+1,s)$.
  Now for each $i$ let $S_i=\{x\in \bar T:x>\bar z^i\}$ so that $S_i\simeq
  T(\alpha,s)$ and every \sss\ of $S$ satisfies \Q.  By PL(3)
  $(u_1,\dots,u_l, w'_i)$ has \P[\e+\delta/2] and we may apply the
  induction hypothesis to obtain a subtree $S'_i$ of $S_i$ isomorphic to
  $T(\alpha,s)$ and having property \P[\e+\delta]\ for $(u_1,\dots,u_l,
  w'_i)$.  It is easy to see that the tree $S=\cup_i(S'_i\cup\{\bar
  z^i\})$ is the required subtree of $T$ isomorphic to $T(\alpha+1,s)$
  with property \P[\epsilon+\delta] for $(u_i)_1^l$.
  
  Let $\alpha$ be a limit ordinal and suppose the result is true for every
  ordinal $\beta<\alpha$.  Let $(\alpha_n)_1^{\infty}$ be the sequence of
  ordinals increasing to $\alpha$ so that
  $T(\alpha,s)=\cup_nT(\alpha_n,s)$.  Let $\epsilon>0, \delta>0$, let
  $(u_i)_1^l$ have \P\ and let $T$ be a tree isomorphic to $T(\alpha,s)$
  satisfying the requirements of the lemma.  Let $S_n$ be the subtree of
  $T$ isomorphic to $T(\alpha_n,s)$ and let $\bar S_n$ be the subtree of
  $S_n$ isomorphic to $T(\alpha_n,s)$ with property 
  \P[\epsilon+\delta/2]\ 
  for $(u_i)_1^l$.  Let $(z^{n,i})_{i=1}^{\infty}$ be the sequence of
  initial nodes of $\bar S_n$ and let $z^{n,i}=(y_{n,i})$ for $n, i\geq1$.
  We have that $(y_{n,i})_{i=1}^{\infty}$ has property
  \P[\epsilon+\delta/2]\ for $(u_i)_1^l$, for each $n\geq1$, and so by
  condition PL(2) we can find sequences $(y_i)_1^{\infty}\subseteq\{
  y_{n,i} : n,i\geq1 \}$ and $1\leq k_1<k_2<\cdots$ with
  $y_i=y_{\phi(i),k_i}$ and such that $(u_1,\dots,u_l,y_1,y_2,\dots)$ has
  property \P[\epsilon+\delta].  Let $S'_n=\{ x\in\bar S_n : x\geq
  (y_{n,k_j})=z^{n,k_j} \mbox{ for some } j\in\phi^{-1}(n) \}$.  Now
  $S'_n$ has \P[\epsilon+\delta]\ for $(u_i)_1^l$ and is still isomorphic
  to $T(\alpha_n,s)$.  We now set $S=\cup_nS'_n$, then $S\simeq
  T(\alpha,s)$ and has property \P[\epsilon+\delta]\ for $(u_i)_1^l$.
\end{proof}

\renewcommand{\labelenumi}{(\roman{enumi})}
\begin{rem}\label{rem:pruning} 
  This is a purely combinatorial result; it could be restated for any set
  $X$.
  \begin{enumerate}
  \item Note that from the construction of $S$ in the proof we have that
    if $z\in S$, then $\{ x\in T : x\leq z \}\subset S$.  In other words,
    if we remove a node $y$ from $T$, then we also remove every node $x\in
    T$ with $x>y$.
  \item If \P\ and \Q\ are such that given $(u_i)_1^l$ with \P\ and
    $(x_i)_1^{\infty}$ with \Q, we can find a subsequence
    $(x'_i)_1^{\infty}$ of $(x_i)_1^{\infty}$ which has property \P\ for
    $(u_i)_1^l$, then we may modify the proof above to remove the
    $\delta$.  Thus, if $T$ has \Q, then we may prune it to obtain $S$
    with \P.
  \item Similarly, if given \Q, we can get \P[\e-\delta] for any
    $\delta>0$, then we may modify the above proof so that given $T$ with
    \Q, we can prune it to obtain $S$ with \P[\e-\delta] for any
    $\delta>0$.
  \end{enumerate}
\end{rem}

\section{Local indices}
\label{sec:l1indices}

In this section we introduce the local indices on a Banach space $X$ that
we shall use throughout this paper.  They have very similar definitions:
one forms trees on $X$ whose nodes satisfy some property $P$, and then the
index is the supremum over the order of the trees.  There are several
different properties that we shall use to produce the different indices.
We first give general results on indices defined in this way, and then we
discuss the specific indices we use.

In the following $X$ will always be a separable Banach space.  Let
$B_{X}=\{ x\in X : \|x\|\leq1 \}$ and $S_{X}=\{ x\in X : \|x\|=1 \}$
denote the unit ball and unit sphere of $X$, respectively.  If $(x_i)_{i\in
  I}$ is a sequence in $X$ for some $I\subseteq\N$, then let $[x_i]_{i\in
  I}$ be the closed linear span of these vectors.  If $X$ also has a basis
$(e_i)_1^{\infty}$, then we define the support of $x\in X$ with respect to
$(e_i)_1^{\infty}$ to be $F\subseteq\N$ if $x=\sum_Fa_ie_i$ with
$a_i\neq0$ for $i\in F$.  If $z=(x_1,\dots,x_n)$ is a sequence of vectors
then $\spp(z)=\cup_1^n\spp(x_i)$.  A sequence $(x_i)_1^{\infty}$ is an
\emph{$\e$ perturbation of a normalized block basis of $(e_i)_1^{\infty}$}
if there exists a normalized block basis $(b_i)_1^{\infty}$ of
$(e_i)_1^{\infty}$ such that $\sum_i\|b_i-x_i\|<\e$.

\begin{defn}\label{defn:P-index}
  Each index will be defined via a property $P$ as follows.  Let $K\geq1$
  and let $P(K)$ be a property, which depends on $K$, that a tree $T$ on
  $X$ may satisfy.  In fact we consider $P(K)$ to be a set of sequences and
  we say that \emph{$T$ is a tree with property $P$ on $X$}, or simply a
  \emph{$P$-tree}, if $T$ is a tree on $X$ with property $P(K)$ for some
  $K\geq1$, i.e.\ for every $(x_i)_1^n\in T$ we have $(x_i)_1^n\in P(K)$.

  For each $K\geq1$, set the \emph{$P(K)$ index} of $X$ to be
  \[ I_P(X,K) = \sup \{o(T) : T \text{ is a tree on $X$ with property }
  P(K) \} \] 
  and then the \emph{$P$ index} of $X$ is given by
  \[ I_P(X) = \sup_{K\geq1} I_P(X,K) \ . \]
\end{defn}

The next theorem contains the general result that if the property $P$ is
sufficiently well behaved, then when the index is countable it will have
the value \wa\ for some $\a<\w_1$.

\begin{thm}\label{thm:property-P}
  Let $K\geq1$, and let $P$ be a property for finite sequences in a Banach
  space satisfying:
  \begin{enumerate}
  \item For every $\xim\in P(K)$, \xim\ is normalized and $K$-basic.
  \item Given $L, C \geq1$ there exists $K'=K'(K,L,C)\geq1$ such that if
    $\xim\in P(K), \yjn \in P(L)$ and $\max(\|x\|,\|y\|) \leq C\|x+y\|$
    for every $x\in[x_i]_1^m, y\in[y_j]_1^n$, then $\xmyn\in P(K')$.
  \item There exists $L=L(K)\geq1$ such that for every
    $\xim\in P(K)$ and any $1\leq k\leq l\leq m$, $(x_i)_k^l\in P(L)$.
  \item There exists $K''=K''(K)\geq1$ such that the closure of $X^n\cap
    P(K)$ in the product topology on $X^n$ is contained in $X^n\cap
    P(K'')$ for every $n\geq1$.
  \end{enumerate}
  Then either $I_P(X)=\w_1$ and there exists $\seq{x_i}\subset S_X$ and
  $K\geq1$ such that $\xim\in P(K)$ for every $m\geq1$, or else
  $I_P(X)=\w^{\a}$ for some $\alpha<\omega_1$.  
\end{thm}

The idea behind the proof is that if we have a $P$-tree on $X$ of order
\ab{\wa}{r} for some $\a<\w_1$ and $r\geq1$, then we can extend this to a
$P$-tree of order \ab{\wa}{(r+1)}.  We do this in
Lemma~\ref{lem:bigger-P-tree} by extending each terminal node of the tree
of order \ab{\wa}{r} with a tree of order \wa.  In order to do this we
show how to concatenate two sequences with property $P$ in the next lemma,
and in Lemma~\ref{lem:delta-tree} we show how to choose a tree of order
\wa\ so that we can use it to extend a finite sequence with property $P$.
Putting all this together gives us the proof of the theorem.

\begin{lem}\label{lem:concat-P(K)-seqs}
  Let $X\subseteq C[0,1]$, let \seq{e_i}\ be a monotone basis for $C[0,1]$
  with basis projections $(P_k)_1^{\infty}$, let $K\geq1$ and let $\xim,
  \yjn\subset X$ be normalized $K$-basic sequences.  If $k\geq 1$ and
  $\delta,\e > 0$ are such that $2\e<\delta$, $\|(I-P_k)x\|\leq\e\|x\|$
  for each $x\in [x_i]^m_1$ and $\|(I-P_k)y\|\geq \delta \|y\|$ for each
  $y\in [y_j]^n_1$, then for every $x\in[x_i]_1^m$ and $y\in[y_j]_1^n$
  \begin{equation}
    \max(\|x\|,\|y\|) \leq \frac{4}{\delta-2\varepsilon}\|x+y\| \
    .\label{eq:star} 
  \end{equation}
\end{lem}

\begin{proof}
  Let $z = \sum^m_1 a_ix_i + \sum^n_1 b_jy_j$, and consider the following
  two possibilities:
  \begin{enumerate}
  \item $\|\sum^n_1 b_jy_j\| \geq \frc 2\|\sum^m_1 a_ix_i\|$;
  \item $\|\sum^n_1 b_jy_j\| < \frc 2\|\sum^m_1 a_ix_i\|$.
  \end{enumerate}
  In case (i),
  \begin{align*}
    \delta \Bigl\|\sum^n_1 b_jy_j\Bigr\| 
           &\leq \Bigl\| (I-P_k)\Bigl( \sum^n_i b_jy_j \Bigr)\Bigr\| \\
           &\leq \|(I-P_k)z\| + \Bigl\| (I-P_k) \Bigl( \sum^m_1 a_ix_i 
                                                       \Bigr) \Bigr\| \\
           &\leq 2\|z\| + \e \Bigl\| \sum^m_1 a_ix_i \Bigr\| \\
           &\leq 2\|z\| + 2\e \Bigl\| \sum^n_1 b_jy_j \Bigr\| \ ,
  \end{align*}
  so that
  \begin{equation}
    \|z\| \geq \frac{\delta-2\e}{2}\Bigl\| \sum^n_1 b_jy_j \Bigr\|
    \geq \frac{\delta-2\e}{4}\Bigl\| \sum^m_1 a_ix_i \Bigr\|
    \label{eq:one} \ .
  \end{equation}
  
  For case (ii)
  \begin{equation}
    \|z\| 
    \geq \Bigl\| \sum^m_1 a_ix_i \Bigr\| - \Bigl\| \sum^n_1 b_jy_j \Bigr\|
    \geq \frc 2 \Bigr\| \sum^m_1 a_ix_i \Bigr\| 
    \geq \Bigl\| \sum^n_1 b_jy_j \Bigr\| \label{eq:two}
  \end{equation}
  and inequality~(\ref{eq:star}) now follows
  {} from inequalities (\ref{eq:one}) and (\ref{eq:two}).

\end{proof}

\begin{lem}\label{lem:delta-tree}
  Let $X\subseteq C[0,1]$ and let \seq{e_i}\ be a monotone basis for
  $C[0,1]$, with basis projections $(P_k)_1^{\infty}$.  Let $K\geq1$ and
  let $T$ be a tree on $X$ of order \wa\ for some $\a<\w_1$ such that each
  node $\xim\in T$ is a normalized $K$-basic sequence.  Then for any
  $k\geq1$ and $0<\delta<1/(2K)$ there exists a tree $T'$ of order \wa\ 
  such that for each node $\yin\in T'$, $\|(I-P_k)y\|\geq \delta \|y\|$
  whenever $y\in [y_i]^n_1$, and for each terminal node $(w_i)_1^n$ of
  $T'$ there exist $l, m\geq1$ and $\xim\in T$ such that $1\leq l \leq
  l+n-1\leq m$ and $w_i = x_{l+i-1}$ for $i=1,\dots,n$.
\end{lem} 
\begin{proof}
  We may assume that $\wa=\lim_n \ab{\w^{\a_n}}{n}$ for some sequence
  $\a_n\nearrow\a$.  We may also assume that $T = \cup_n T(n)$ with $T(n)$
  isomorphic to the replacement tree $T(n,\w^{\a_n})$.  It is sufficient
  to find a sequence $n_r\nearrow\w$ and trees $T'_r\subseteq T(n_r )$ of
  order $\ab{\w^{\a_{n_r}}}{r}$ such that $T'_r$ satisfies the conditions
  of the lemma for each $r\geq1$.
    
  Let $0 < \xi < 1/(2K) - \delta$ and let $N \geq1$ and sets $(A_l)^N_1$
  satisfy for $1\leq l\leq N$:
  \[ A_l\subseteq B_{[e_i]^k_1},\quad \diam A_l < \xi,\quad 
     \cup^N_1 A_l = B_{[e_i]^k_1} \ . 
  \]
  Let $r\geq 1$ and choose $n\geq (N + 1)r$.  Consider a subtree $S_1$ of
  $T(n)$ isomorphic to $T(N + 1, \ab{\w^{\a_n}}{r} )$ and let $F_1 :
  R(S_1)\rightarrow T_{N+1} = \{a_1,\dots, a_{N+1} \}$, with 
  $a_1 <\dots< a_{N+1}$, be the defining map for the replacement tree.
  Then $F_1^{-1}(a_1 )$ is isomorphic to $T_{\w^{\a_n}\cdot r}$.  If $\|(I
  - P_k )y\| \geq\delta\|y\|$ for all $y \in [y_i]^j_1$, and every
  $(y_i)^j_1\in F_1^{-1}(a_1 )$, then $F_1^{-1} (a_1 )$ is the subtree we
  seek.  If not, then there exists a terminal node $(y^1_i)^{k_1}_1$ of
  $F_1^{-1}(a_1 )$ and $y_1 \in [y^1_i]^{k_1}_1$ with $\|y_1 \| = 1$ and
  $\|(I - P_k )y_1 \| < \delta$.
  
  Let $S_2 = \{(u_1 ,\dots, u_j) \in S_1 : j >
  k_1$ and $u_i = y^1_i\ (i = 1,\dots, k_1 )\}$, so that the restricted
  tree $R(S_2)$ is isomorphic to $T(N, \ab{\w^{\a_n}}{r} )$.  Let $F_2 :
  R(S_2) \rightarrow T_N = \{a_2 ,\dots, a_{N+1} \}$ be the restriction
  of $F_1$ to $R(S_2)$ so that $F_2^{-1}(a_2 ) \simeq T_{\w^{\a_n}\cdot
    r} $.  Again, either $F_2^{-1}(a_2 )$ is the required tree, or else
  there is a terminal node $(y^2_i)^{k_2}_1$ in $F_2^{-1}(a_2 )$ and $y_2
  \in [y^2_i]^{k_2}_1$ with $\|y_2 \| = 1$ and $\|(I - P_k )y_2 \| <
  \delta$.  Continuing in this way we obtain either a subtree $T'_r$
  isomorphic to $T_{\w^{\a_n}\cdot r}$ such that $\|(I - P_k )y\| \geq
  \delta\|y\|$ for each $y \in [y_i]^j_1$, and every $(y_i)^j_1 \in T'_r$,
  or else there is a branch $(y_1^1,\dots, y^1_{k_1},\dots, y^{N+1}_1
  ,\dots, y^{N+1}_{k_{N+1}})$ of $T$ and normalized vectors $y_j \in
  [y^j_i]_{i=1}^{k_j}$ such that $\|(I - P_k )y_j\| < \delta$ for $j =
  1,\dots, N + 1$.  But then there exist 
  $l, j, j' \in \{1,\dots, N + 1\}\ 
  (j \neq j')$ such that $P_k y_j, P_k y_{j'}\in A_l$, and hence
  \[ \frc K \leq \|y_j - y_{j'}\| \leq  \|P_k (y_j - y_{j'})\| +
                 \|(I - P_k )(y_j - y_{j'})\| 
            < \xi + 2\delta < \frc K \ , 
  \]
  a contradiction.  The last condition on $T$ is clear from the proof.
\end{proof}

\begin{lem}\label{lem:bigger-P-tree}
  Let $X\subseteq C[0,1]$ and let $P$ be a property satisfying conditions
  (i)--(iv) of Theorem~\ref{thm:property-P}.  For all $K\geq1$ there
  exists $L\geq1$ such that for evey $\alpha<\omega_1$ and $r\geq1$, if
  there exists a $P(K)$ tree on $X$ of order \ab{\wa}{r}, then there
  exists a $P(L)$ tree on $X$ of order \ab{\wa}{(r+1)}.
\end{lem}
\begin{proof}
  Let $S$ be the given $P(K)$ tree on $X$ of order \ab{\wa}{r}\ and let
  $T$ be a $P(K)$ tree on $X$ of order \wa{}.  Let $(z^j)^{\infty}_1$ be
  the sequence of terminal nodes of $S$, so that $z^j =
  (x^j_i)^{m_j}_{i=1}$.  Choose $0<\delta<1/(2K)$, $0<\e<\delta/2$ and for
  each $j\geq1$ find $k_j\geq1$ such that $\|(I-P_{k_j})x\|\leq\e\|x\|$
  whenever $x\in[x^j_i]_1^{m_j}$.  Apply the previous lemma to $T$ for $K$
  and $\delta$ with $k=k_j$ to obtain a tree $T(z^j)$ of order
  \wa{} such that for each node $\yin\in T(z^j)$ we have
  $\|(I-P_{k_j})y\|\geq\delta\|y\|$ whenever $y\in[y_i]_1^n$.  From
  condition (iii) of Theorem~\ref{thm:property-P} and the construction of
  $T(z^j)$ there exists $K'=K'(K)$ such that $T(z^j)$ is a $P(K')$ tree
  and hence by Lemma~\ref{lem:concat-P(K)-seqs} and condition (ii) for
  property $P$, there exists $L=L(K,\delta,\varepsilon)$ such that the
  tree
  \[ S(z^j) = 
    \{(x^j_1,\dots,x^j_{m_j},y_1,\dots,y_n) : (y_i)^n_1 \in T(z^j)\} \cup 
    \{(x^j_1),\dots,(x^j_1,\dots,x^j_{m_j})\}
  \]  
  is a $P(L)$ tree on $X$.  The tree $\cup_{j=1}^{\infty}S(z^j)$ is the
  required $P(L)$ tree on $X$ of order \ab{\wa}{(k + 1)}.
\end{proof}

\begin{proof}[Proof of Theorem~\ref{thm:property-P}]
  If $I_P(X)=\w_1$, then, since the closure of a $P(K)$ tree is a $P(K')$
  tree for some $K'\geq1$ by condition (iv), it follows from
  Proposition~\ref{prop:wf=>cble} that there exists an infinite sequence
  as in the statement of the theorem.
  
  Otherwise we assume the index is countable and let T be a $P(K)$-tree on
  $X$ of order $\w^{\g}$.  By the previous lemma there exist numbers
  $K_i\geq 1$ and $P(K_i)$-trees $T_i$ on $X$ of order \ab{\w^{\g}}{i}\ 
  for $i = 1, 2,\dots\ (K_1 = K)$.  Therefore the $P$-index is at least
  $\w^{\g+1}$.  It follows that the $P$-index is
  \begin{align*}
    I_P(X) &= \sup \{ \ab{\w^{\g}}{k} : \text{ there exists } K \text{ and
    a $P(K)$-tree on $X$ of order } \w^{\g}, k \in \N\} \\
      &= \sup \{ \w^{\g+1} : \text{ there exists $K$ and a $P(K)$-tree on
    X of order } \w^{\g}\} \\  
      &= \wa{}
  \end{align*}  
for some $\a < \w_1$.
\end{proof}

\begin{defn}\label{defn:properties-P}
  We shall use the following indices; we give the name of the index, the
  symbol we use for it and then the property $P$ that each node of the
  tree must satisfy for $K\geq1$.
  \begin{entry}
    \item [\lli, $I(X)$] \xim\ is normalized and $K$-equivalent to the
      unit vector basis of $\ell_1^m$, $\xim\kequiv\ell_1^m$, i.e.
      \[ \frac{1}{K}\sum_1^m|a_i| \leq \Bigl\|\sum_1^ma_ix_i\Bigr\| \leq
      \sum_1^m|a_i| \] 
      for every $(a_i)_1^m\subset {\mathbf R}$.  (Of course, this second
      inequality is always true.)
    \item [\llpi, $I^+(X)$] \xim\ is an \emph{\llpk-sequence}, i.e. \xim\ 
      is normalized, $K$-basic and satisfies
      \[  \frac{1}{K}\sum_1^m a_i \leq \Bigl\| \sum_1^ma_ix_i \Bigr\| \]
      whenever $(a_i)_1^m\subset\R^+$.
    \item [\lii, $J(X)$] \xim\ is normalized and $K$-equivalent to the
      unit vector basis of $\ell_{\infty}^m$,
      $\xim\kequiv\ell_{\infty}^m$, i.e.  there exist $c,C\geq1$ such that
      $cC\leq K$ and
      \[ \frc{c} \max_{1\leq i\leq m}|a_i| \leq \Bigl\| \sum_1^ma_ix_i
      \Bigr\| \leq C \max_{1\leq i\leq m}|a_i| \]
      for every $(a_i)_1^m\subset {\mathbf R}$.
    \item [\lipi, $J^+(X)$] \xim\ is an \emph{\lipk-sequence}, i.e. \xim\ 
      is normalized, $K$-basic and there exists a sequence
      $(a_i)_1^m\subset[1/K,1]$ such that $\|\sum_1^ma_ix_i\|\leq K$.
  \end{entry}
\end{defn}

\begin{rem}\label{rem:l1p-equiv}
  It is an easy consequence of the geometric Hahn-Banach theorem that the
  following is an equivalent definition of an \llp sequence.  Let
  $K\geq1$, then a sequence $(x_i)\subset X$ (finite or infinite) is an
  \emph{\llpk-sequence} if and only if $(x_i)$ is a normalized, $K$-basic
  sequence and there exists $f\in S_{X^*}$ such that
  $f(x_i)\geq\frac{1}{K}$ for each $i$.
\end{rem}

Rosenthal~\cite{r2} studied other aspects of \llp\ and \lip\ sequences
under the names wide-($s$) and wide-($c$) sequences respectively, with
different quantifications (see Lemma~\ref{lem:plus<->wide}).  For
$\lambda\geq1$, \xim\ is a \emph{$\lambda$-wide-($s$)} sequence in a
Banach space $X$ if it is $2\lambda$-basic with $\|x_i\|\leq \lambda$ for
$i\leq m$ and
\[ \sup_{k\leq m} \Bigl| \sum_k^m a_i \Bigr| \leq \lambda
   \Bigl\| \sum_1^ma_ix_i \Bigr\| 
\]
for every $(a_i)_1^m\subset {\mathbf R}$.  A \emph{$\lambda$-wide-($c$)}
sequence in $X$ is a $\lambda$-basic sequence $\xim\subset X$ with
$1/\lambda \leq \|x_i\|\leq 1$ for $i\leq m$ and $\| \sum_1^m x_i \|
\leq \lambda$.  We investigate the relationships between these notions in
Lemma~\ref{lem:plus<->wide} below.

Each of the above indices has a companion \emph{block basis index}, where
the property has the additional requirement that the space $X$ have a
basis \seq{e_i} and that each node \xim\ of the tree be a block basis of
\seq{e_i}.  The block basis indices are written $I_b(X)$, $I^+_b(X)$ etc.
They are calculated with respect to a fixed basis $(e_i)_1^{\infty}$ of
$X$ and should more properly be written $I_b(X,(e_i)_1^{\infty})$ and
$I^+_b(X,(e_i)_1^{\infty})$ etc., since they depend on the basis.  For
example James' space, $J$, has a shrinking basis $(e_i)_1^{\infty}$ and
$I_b(J,(e_i)_1^{\infty})<\w_1$, but it also has a non-shrinking basis
$(f_i)_1^{\infty}$, and for this we have $I_b(J,(f_i)_1^{\infty})=\w_1$.
For the \lli, from \cite{jo} we know that if $I(X)=\w^{1+\a}$ for some
$\a<\w_1$, then $I_b(X)=\w^{\a}$ or $\w^{1+\a}$.  It is easy to construct
spaces X (see~\cite{jo} Remark~5.15 (ii)) with two different bases
$(e_i)_1^{\infty}$ and $(f_i)_1^{\infty}$, such that $I(X)=\w^{n+1}$,
$I^+_b(X,(e_i)_1^{\infty})=\w^n$ and $I_b(X,(f_i)_1^{\infty})=\w^{n+1}$.
However, because we shall be working with a fixed basis for $X$ we shall
omit reference to the basis.  Furthermore, if the basis for $X$ is
unconditional, then $I^+_b(X)=I_b(X)$.

The trees used to calculate each index are named after the index.  Thus a
tree with property $P(K)$ for the \lli\ is called an \llkt, or just an
\llt.

\begin{cor}[Corollary to
  Theorem~\ref{thm:property-P}]\label{cor:I-P-bb-index} 
  Let $P$ be a property satisfying the conditions of
  Theorem~\ref{thm:property-P}.  Let \seq{e_i} be a basis for a Banach
  space $X$, and let $I_P(X,\seq{e_i})$ be a block basis index defined via
  $P$.  Either $I_P(X,\seq{e_i})=\w_1$, and there exist $K\geq1$ and a
  normalized block basis \seq{x_i} of \seq{e_i} such that $\xim\in P(K)$
  for every $m\geq1$, or else $I_P(X,\seq{e_i})=\w^{\a}$ for some
  $\a<\w_1$.  
\end{cor}
\begin{proof}
  We only need make a couple of modifications to the proof of
  Theorem~\ref{thm:property-P}.  Instead of embedding $X$ into $C[0,1]$
  and using a basis there we use the basis \seq{e_i} of $X$.  Then in the
  proof of Lemma~\ref{lem:bigger-P-tree} we must ensure that we construct
  a block basis tree.  But this is easy.  In
  Lemma~\ref{lem:bigger-P-tree}, for each terminal node $z^j =
  (x^j_i)^{m_j}_{i=1}$ we choose $k_j$ so that $P_{k_j}x^j_i=x^j_i\ 
  (i=1,\dots,m_j)$, and once we find the subtree $T(z^j)$ of $T$, we take
  a further subtree
  \[ T(z^j)' = \{ (x_i)_{k_j+1}^m : k_j<m, \xim\in T(z^j) \} \ , \]
  which also has order \wa.  Since we started with a block basis tree $T$,
  the subtrees $T(z^j)$ and $T(z^j)'$ will also be block basis trees, as
  will the final tree.
\end{proof}

We next note for future reference that each of the four properties we have
defined are easily seen to satisfy the conditions of
Theorem~\ref{thm:property-P}.  The proofs are elementary calculations.

\begin{lem}\label{lem:properties-good}
  Each of the four properties in Definition~\ref{defn:properties-P}
  satisfies conditions (i)--(iv) in Theorem~\ref{thm:property-P} for a
  property $P$.  
\end{lem}

We have the following Corollary of Theorem~\ref{thm:property-P}, which
includes some results already proved in~\cite{jo}:

\begin{cor}\label{cor:indices-all-wa}
  For a separable Banach space $X$, each of the indices $I(X)$, $I_b(X)$, 
  $I^+(X)$, $I^+_b(X)$, $J(X)$, $J_b(X)$, $J^+(X)$, $J^+_b(X)$ is either
  uncountable or \wa\ for some $\alpha<\omega_1$.
\end{cor}

\begin{defn}\label{defn:i-plus-w}
  We need one more index which doesn't fit into this pattern since it
  relies on the structure of the trees used.  For $K\geq1$ a tree $T$ on
  $X$ is an \llpkwt\ if each node $\xim\in T$ is an \llpks\ and $T$ is a
  weakly null tree (see Definition~\ref{defn:wnt}).  The \llpwi\ is
  written $I^+_w(X)$.  We show in Theorem~\ref{thm:I+wX=w^a} that when
  $I^+_w(X)$ is countable it is also equal to \wa\ for some $\a<\w_1$.
\end{defn}

In the light of the block basis indices above we make the following
definition.  

\begin{defn}\label{defn:bbt}
  Let $X$ be a Banach space with a basis $(e_i)_1^{\infty}$.  A
  \emph{block basis tree on $X$} is a tree on the unit sphere $S_{X}$ of
  $X$ such that every node $(x_i)_1^n\in T$ is a block basis of
  $(e_i)_1^{\infty}$.
\end{defn}

As well as taking subtrees we also want to take ``block trees''.  They are
an extension of the notion of a block basis to trees.

\begin{defn}\label{defn:blockt}
Let $T$ be a tree on the unit sphere $S(X)$ of a Banach space $X$.  We
say $S$ is a \emph{block tree} of $T$, written $S\preceq T$, if $S$
is a tree on $S(X)$ such that there exists a subtree $T'\subset T$ and
an isomorphism $f: T'\stackrel{\sim}{\rightarrow}S$ satisfying:
$f((x_i)_1^m)=(y_i)_1^n$ is a normalized block basis of $(x_i)_1^m$
for each $(x_i)_1^m\in T'$, and if $(x_i)_1^k\in T'$ for some $k<m$,
then there exists $l<n$ such that $(y_i)_1^l=f((x_i)_1^k)$ and
$(y_i)_{l+1}^n$ is a normalized block basis of $(x_i)_{k+1}^m$.
\end{defn}

In the next theorem we give some basic properties of these indices.
Statement (\ref{thm:j-reflexive:item:A}) is due to Bourgain~\cite{b}, and
(\ref{thm:j-reflexive:item:B}), (\ref{thm:j-reflexive:item:C}) were proven
in~\cite{jo}.
\begin{thm}\label{thm:j-reflexive}
  Let $X$ be a separable Banach space, then
  \begin{enumerate}[\rm (i)]
  \item\label{thm:j-reflexive:item:A} $I(X)<\w_1$ if and only if $\ell_1$
    does not embed into $X$;
  \item\label{thm:j-reflexive:item:B} If $I(X)\geq\w^{\w}$, then
    $I(X)=I_b(X)$;
  \item\label{thm:j-reflexive:item:C} If $I(X)=\w^n$ for some $n\in\N$,
    then $\ib{X}=\w^m$ where $m\in\{n,n-1\}$;
  \item\label{thm:j-reflexive:item:D} $I^+(X)<\w_1$ if and only if $X$ is
    reflexive;
  \item\label{thm:j-reflexive:item:E} $I^+(X)=\w$ if and only if $X$ is
    super-reflexive;
  \item\label{thm:j-reflexive:item:F} $I^+(X)=J^+(X)$ and hence
    $J^+(X)<\w_1$ if and only if $X$ is reflexive;
  \item\label{thm:j-reflexive:item:G} $I^+_b(X,(e_i)_1^{\infty})<\omega_1$
    if and only if $(e_i)_1^{\infty}$ is a shrinking basis for $X$;
  \item\label{thm:j-reflexive:item:H} $J_b^+(X,\seq{e_i})<\w_1$ if and
    only if \seq{e_i} is a boundedly complete basis for $X$;
  \item\label{thm:j-reflexive:item:I} If $\ell_1\not\hookrightarrow X$,
    then $I^+_w(X)<\omega_1$ if and only if $X^*$ is separable.
  \end{enumerate}
  $X$ is assumed to have a basis $(e_i)_1^{\infty}$ in
  (\ref{thm:j-reflexive:item:B}), (\ref{thm:j-reflexive:item:C}),
  (\ref{thm:j-reflexive:item:G}) and (\ref{thm:j-reflexive:item:H}).
\end{thm}

\begin{proof}
  Statements (\ref{thm:j-reflexive:item:D})
  and (\ref{thm:j-reflexive:item:E}) are results of James~\cite{j2} and
  Milman and Milman~\cite{mm} stated in terms of the \llpi.  In particular
  (\ref{thm:j-reflexive:item:D}) follows from \cite{j2}~Theorem~1 or
  \cite{mm}~Corollary of Theorem~2, while (\ref{thm:j-reflexive:item:E})
  follows from \cite{j3} and \cite{j4}.
  
  Before we can give the proof of (\ref{thm:j-reflexive:item:F}) we need
  some results on the relationship between \llp\ and \lip\ sequences, so
  we shall postpone the proof until we have these..
  
  For part (\ref{thm:j-reflexive:item:G}), if $(e_i)_1^{\infty}$ is not
  shrinking, then it has a normalized block basis which is an \llp\ 
  sequence, giving an \llpbbt\ of order $\w_1$, and the other direction
  follows from Proposition~\ref{prop:wf=>cble}, since $X$ is separable and
  the closure of an \llpbbt\ is again an \llpbbt.
  
  If $J_b^+(X,\seq{e_i})=\w_1$, then there exist $K\geq1$ and a normalized
  block basis \seq{x_i} of \seq{e_i} so that \xim\ is an \lipk\ sequence
  for each $m\geq1$.  By a compactness argument it is easy to find
  $\seq{a_i}\subset[1/K,1]$ such that $\|\sum_1^m a_ix_i \|\leq 2CK$,
  where $C\geq1$ is the basis constant of \seq{e_i}.  Thus \seq{e_i} is
  not boundedly complete.  The converse is clear and part
  (\ref{thm:j-reflexive:item:H}) follows.
  
  The proof of part (\ref{thm:j-reflexive:item:I}) requires more work.
  First, if $X^*$ is separable, then by Zippin,~\cite{z} $X$ embeds into a
  Banach space with a shrinking basis.  Thus it is sufficient to prove
  that if $X$ has a shrinking basis $(e_i)_1^{\infty}$, then
  $I^+_w(X)<\w_1$.  If we show that $I^+_w(X)\leq I^+_b(X)$, then this
  would follow from part (\ref{thm:j-reflexive:item:G}).  To show this we
  will take an \llpwt\ on $X$ of order \a, and apply the Pruning Lemma to
  obtain a perturbation of an \llpbbt\ on $X$ of order \a.  Let $T$ be an
  \llpwt\ on $X$ of order \a, so that $T$ is isomorphic to $T(\a,s)$.
  Define $(x_i)_1^{\infty}$ to have property \Q\ if it is weakly null, and
  define $(x_i)_1^{\infty}$ to have property \P\ if it is an
  $\e$-perturbation of a block basis of $(e_i)_1^{\infty}$.  Using that
  $(e_i)_1^{\infty}$ is shrinking it is standard work to show that \Q\ and
  \P\ satisfy the requirements of the Pruning Lemma, and hence we may
  obtain a subtree $T'$ of order \a\ which is a perturbation of an
  \llpbbt\ of order \a\ on $X$, i.e. each terminal node is an
  \e-perturbation of an \llp\ block basis of $(e_i)_1^{\infty}$.
  
  Next suppose that $X^*$ is not separable; we shall show that
  $I^+_w(X)=\w_1$.  Let $\Delta$ denote the Cantor set, and let
  $(A_{n,i})$ be a sequence of subsets of $\Delta$ for $n=0,1,2,\dots$ and
  $i=0,1,2,\dots,2^n-1$ such that $A_{0,0}=\Delta$ and each $A_{n,i}$ is
  the union of the disjoint, non-empty, clopen sets $A_{n+1,2i}$ and
  $A_{n+1,2i+1}$, with $\lim_{n\rightarrow\infty} \sup_{0\leq i<2^n}
  \diam(A_{n,i}) = 0$.  A \textit{Haar system on $\Delta$} (relative to
  $(A_{n,i})$) is a sequence of continuous functions,
  $(h_m)_1^{\infty}\subseteq C(\Delta)$ with
  \[ h_{2^n+i} = \ch_{A_{n+1,2i}}-\ch_{A_{n+1,2i+1}} (n=0,1,2,\dots,\ 
  i=0,1,2,\dots,2^n-1) \ , \] 
  where $\ch_A$ is the characteristic function of the set $A$.  A sequence
  of continuous functions, $(g_m)_1^{\infty}\subseteq C(\Delta)$ is a
  \textit{Haar system on $\Delta$ up to $\zeta$} if there exists a Haar
  system $(h_m)_1^{\infty}$ on $\Delta$ such that for each $m\geq1$, $\spp
  g_m=\spp h_m$, $\sign g_m = \sign h_m$ and $1-\zeta\leq|g_m(x)|\leq1$
  for every $x\in\spp g_m$.
  
  Given a Banach space $X$, and a subset $\Delta\subseteq B_{X^*}$ which
  is weak$^*$ homeomorphic to the Cantor set, a sequence
  $(x_m)_1^{\infty}\subseteq S_X$ is a \textit{Haar system up to $\zeta$,
    relative to $\Delta$} if the restrictions $(x_m|_{\Delta})_1^{\infty}$
  form a Haar system up to $\zeta$ on $C(\Delta)$.
  
  By a result of Stegall~\cite{st}, since $X^*$ is not separable, we have
  that given $\zeta>0$ there exists a set $\Delta\subseteq B_{X^*}$ which
  is weak$^*$ homeomorphic to the Cantor set, and a dyadic tree $S$ of
  elements in $S_X$ which form a Haar system up to $\zeta$, relative to
  $\Delta$.  The dyadic tree is the natural one, with $x<y$ if and only if
  $\spp x|_{\Delta} \supset \spp y|_{\Delta}$.  Let $\tau_{\w} = \{
  (n_1,\dots,n_k) : 1\leq k,\ n_i\in\N,\ i\leq k \}$ be the countably
  branching tree with \w\ levels, ordered by $(m_1,\dots,m_k)\leq
  (n_1,\dots,n_l)$ if and only if $k \leq l$ and $m_i = n_i$ for $i\leq
  k$.  Choose $S_1\subseteq S$ isomorphic to $\tau_{\w}$ so that if $x\in
  S_1$ has immediate successors $(x_i)_1^{\infty}$ (i.e.\ if $x$ is
  equivalent to the node $(n_1,\dots,n_k)$ of $\tau_{\w}$, then $x_i$
  corresponds to $(n_1,\dots,n_k,i)$ for $1\leq i$), then
  \[ \spp x_i|_{\Delta} \subseteq\{ x^*\in\Delta: x(x^*)\geq1-\zeta \} \
  . \] 
      
  Since $X$ does not contain $\ell_1$, it follows by~\cite{r} that we may
  prune $S_1$ to obtain a further subtree $S_2$ isomorphic to $\tau_{\w}$
  so that if $(x_i)_1^{\infty}$ is the sequence of immediate successors of
  $x\in S_2$ as above (or the sequence of initial nodes), then
  $(x_i)_1^{\infty}$ is weak-Cauchy.  Furthermore they will still satisfy
  the property above and we may assume that if we set
  \[ y_j = (x^{2j-1}-x^{2j})/\|x^{2j-1}-x^{2j}\| \  \mbox{ for }j\geq1 \ ,
  \] 
  then $(y_i)_1^{\infty}$ is 2-basic and weakly null.  From this we can
  create a tree $S_3\subseteq S_X$ which is isomorphic to $\tau_{\w}$ with
  nodes $y_i$ as above and letting the immediate successors of such a node
  $y_i$ be formed in the same manner from the successors in $S_2$ of
  $x_{2i-1}$.
  
  The resulting tree has the property that if $(z_i)_1^{\infty}$ is a
  branch of $S_3$, then $\spp z_{i+1}|_{\Delta} \subset \{ x^*\in\Delta :
  z_i(x^*)>(1-\zeta)/2 \}$ for each $i\geq1$.  Hence for each such branch
  there exists $x^*\in B_{X^*}$ such that $x^*(z_i)>(1-\zeta)/2$ for every
  $i\geq1$ so that $(z_i)_1^{\infty}$ is an \llp\ sequence.  Now, for all
  $\a<\w_1$, since $S_3$ is isomorphic to $\tau_{\w}$, it follows that we
  may find a subtree of $S_3$ which is a weakly null tree isomorphic to
  $T(\a,s)$ and so $I^+_w(X)=\w_1$ as claimed.
\end{proof}

\begin{rem}
  The condition that $\llo\not\hookrightarrow X$ is necessary for
  (\ref{thm:j-reflexive:item:I}) above.  Indeed, $I^+_w(\llo)=0$ since
  there are no non-empty weakly null trees on \llo.
\end{rem}

We now present the relationship between Rosenthal's \ws\ and \wc\ 
sequences and our \llp\ and \lip\ sequences, in order to prove part
(\ref{thm:j-reflexive:item:F}) of the above theorem.  We first note the
following result on the relationship between \ws\ and \wc\ sequences:
\begin{prop}[Rosenthal~\cite{r2}]\label{prop:R}
  Let $(b_j)$ be a sequence in $X$, finite or infinite, and let $(e_j)$ be
  its difference sequence, so that $e_1=b_1$ and $e_j=b_j-b_{j-1}\ (j>1)$.
  Then for every $\lambda\geq1$ there exists $\mu\geq1$ such that
  \begin{enumerate}
  \item if $(b_j)$ is \lws, then $(e_j)$ is \mwc;
  \item if $(e_j)$ is \lwc, then $(b_j)$ is \mws.
  \end{enumerate}
\end{prop}

In the next lemma we show how to move between \llp\ and \ws\ sequences,
and between \lip\ and \wc\ sequences, and then in the lemma following we
show how one can perturb this process and still keep control of the
constants.  This will be important when moving between these sequences in
trees.  We leave the proofs as exercises.

\begin{lem}\label{lem:plus<->wide}
  Let $\lambda, K \geq 1$, and let $\xim\subset X$.
  \begin{enumerate}
  \item If \xim\ is an \llpk\ sequence and $f\in S_{X^*}$ is such that
    $f(x_i)\geq1/K$ for $1\leq i \leq m$, then $(x_i/f(x_i))_1^m$ is a
    2K-\wss.  
  \item If \xim\ is a \lws\ sequence, then $(x_i/\|x_i\|)_1^m$ is an
    \llp\  sequence with constant $\max(2\lambda, \lambda^2)$.
  \item If \xim\ is an \lipk\ sequence and $(b_i)_1^m\subset[1/K,1]$ is
    such that $\|\sum_1^m b_ix_i\|\leq K$, then $(b_ix_i)_1^m$ is a
    $K$-\wcs.
  \item If \xim\ is a \lwc\ sequence, then $(x_i/\|x_i\|)_1^m$ is an
    \lip-$\lambda$ sequence.
  \end{enumerate}
\end{lem}

\begin{lem}\label{lem:plus<->wide-upto-e}
  Let $K\geq1$ and $\e>0$.  Then there exist $\e_i\searrow0$ such that for
  $\xim\subset X$,  
  \begin{enumerate}
  \item if \xim\ is an \llpk\ sequence, $f\in S_{X^*}$ satisfies
    $f(x_i)\geq 1/K$ for $i\leq m$, and $(c_i)_1^m\subset[1/K,1]$ is
    chosen so that $|f(x_i)-c_i|<\e_i$ for $i\leq m$, then $(x_i/c_i)_1^m$
    is $2(K+\e)$-\ws;
  \item if \xim\ is an \lipk\ sequence, $(b_i)_1^m\subset[1/K,1]$
    satisfies $\|\sum_1^mb_ix_i\|\leq K$, and $(c_i)_1^m\subset[1/K,1]$ is
    chosen so that $|b_i-c_i|<\e_i$ for $i\leq m$, then $(c_ix_i)_1^m$ is
    $(K+\e)$-\wc.
  \end{enumerate}
\end{lem}

The next lemma gives the framework for applying the previous two results
to swap between whole trees of these sequences.

\begin{lem}\label{lem:approx-trees}
  Let $\a<\w_1$, $\delta>0$ and $\seq{\e_i}\subset\R^+$.  Let $T$ be a
  tree of order \a, \seq{z^i} the sequence of terminal nodes of $T$, with
  $z^i=(x^i_j)_{j=1}^{m_i}$, and for each $i$, let $f^i$ be a map {} from
  $\{x^i_1,\dots,x^i_{m_i}\}$ into $[\delta,1]$.  Then there exists an
  increasing sequence $N \subseteq \N$ such that the subtree $\bar T = \{
  z\in T : z \leq z^{n} \text{ for some } n\in N \}$ of $T$ has order \a\ 
  and $| f^{n}(x^{n}_l) - f^{n'}(x^{n'}_l)| < \e_l$ whenever $n, n'\in N$
  and $x^{n}_k=x^{n'}_k$ for $1\leq k \leq l$.
\end{lem}
\begin{proof}
  We prove this by induction on \a.  There is nothing to prove for
  $\a\leq\w$, and if the result has been proven for every $\b<\a$, then it
  is also clear when \a\ is a limit ordinal.  Thus suppose the result has
  been proven for $\a'$, let $\a=\a'+1$, and let $T$ be a tree of order
  \a.  By taking a subtree we may assume that $T$ has a unique initial
  node, $z=(x_1)$, and let $S=\{(x_j)_2^m:(x_j)_1^m\in T\}$, a tree of
  order $\a'$.  The terminal nodes of $S$ are $w^i=(x^i_j)_{j=2}^{m_i}$,
  and let $N\subseteq\N$ be the sequence from the induction hypothesis on
  $S$ for $\a', \delta$ and the sequence $(\e_i)_{i=2}^{\infty}$.  Now let
  $T'$ be the subtree of $T$, $\{z\in T : z\leq z^{n} \text{ for some }
  n\in N \}$.  We have that $|f^{n}(x^{n}_l) - f^{n'}(x^{n'}_l)| < \e_l$
  whenever $n, n'\in N$, $l\geq2$ and $x^{n}_k=x^{n'}_k$ for
  $2\leq k \leq l$.  We must now stabilize the maps on $x_1=x_1^n$ for
  each $n\in N$.

  Let $m=[(1-\delta)/\e_1]+1$ and for $1\leq r\leq m$ let 
  \[ M_r = \{ n\in N : f^{n}(x_1)\in [\delta+(r-1)\e_1,\delta+r\e_1) \} \
  . \]  
  This forms a partition of the terminal nodes of $T'$ and by
  \cite{jo}~Lemma 5.10 one of the trees
  \[ T'(M_{r}) = \{ y\in T':y\leq z^n \text{ for some } n\in M_{r}\} \]  
  has order \a\ for some $1\leq r_0 \leq m$.  The sequence
  $M_{r_0}\subseteq N \subseteq\N$ is now the required sequence.
\end{proof}

\begin{proof}[Proof of
  Theorem~\ref{thm:j-reflexive}~(\ref{thm:j-reflexive:item:F})]
  We shall show that for each $\a<\w_1$ there exists an \llpt\ on $X$ of
  order \a\ if and only if there exists an \lipt\ on $X$ of order \a.
  
  Let $\a<\w_1$ and let $T$ be an \llpkt\ on $X$ of order \a\ for some
  $K\geq1$.  Choose $0<\e\ll 1/K$ and a sequence $\seq{\e_i}\subset(0,\e)$
  decreasing rapidly to zero.  Let $\bar T$ be the tree obtained when
  Lemma~\ref{lem:approx-trees} is applied to $T$ with $\delta=1/K$ where
  the functions $f^i$ on $T$ are in $S_{X^*}$ with $f^i(x^i_j)\geq1/K\ 
  (1\leq j\leq m_i)$ for each terminal node $(x^i_j)_{j=1}^{m_i}$.  We can
  find these functions since the terminal nodes are \llpk\ sequences.  Let
  $f:\{ x^i_j : 1\leq j\leq m_i, i=1,2,\dots \} \rightarrow [1/K,1]$ be
  the map $f(x^i_j)=f^k(x^i_j)$ where $k=\min\{i'\geq1 :
  x^{i'}_j=x^i_j\}$ and let
  \[ S = \bigcup_{i=1}^{\infty} \Bigl\{
         \Bigl(\frac{x^i_1}{f(x^i_1)}\Bigr),\dots,
         \Bigl(\frac{x^i_1}{f(x^i_1)},\dots,
               \frac{x^i_{m_i}}{f(x^i_{m_i})}\Bigr) \Bigr\} \ ,
  \] 
  Clearly $S$ has order \a\
  and is a $2(K+\e)$-\wst\ by Lemma~\ref{lem:plus<->wide-upto-e}.
  
  Now let $U=\{ (y_1,y_2-y_1,\dots,y_n-y_{n-1}) : \yjn\in S \}$, then $U$
  is a \wc\ tree by Proposition~\ref{prop:R} above and clearly $U$ is
  isomorphic to $S$ so that $o(U)=\a$.  Finally let $V=\{
  (u_i/\|u_i\|)_1^m : (u_i)_1^m \in U \}$, then $V$ is an \lipt\ by
  Lemma~\ref{lem:plus<->wide}~(iv) and isomorphic to $U$ so it is of order
  \a.  This completes the proof of one implication.  The other is similar.
\end{proof}

\begin{rem}\label{rem:wide-s=wide-c}
  Notice that this proof also shows that the \ws\ and \wc\ indices are both
  equal to $I^+(X)$.
\end{rem}

In Theorem~1.1 of~\cite{jo} the following extension of the finite version
of the result of James~\cite{j} that $\ell_1$ is not distortable is shown.
Given $K\geq1$, $\varepsilon>0$ and $\alpha<\omega_1$ there exists
$\beta<\omega_1$ such that if $T$ is an \llkt\ on $X$ of order \b, then
there exists a block tree $T'$ of $T$ which is an \llt\ with constant
$1+\e$ and order \a.  We have the analogous results for \llpbbt s and for
\llpwt s:

\begin{thm}\label{thm:localization}
  Let $K\geq1$, $\varepsilon>0$ and $\alpha<\omega_1$.
  \begin{enumerate}[\rm (i)]
  \item There exists $\beta<\omega_1$ such that for every \llpkwt\ $T$ of
    order \b\ (i.e.\ isomorphic to $T(\b,s)$) on a Banach space with
    separable dual there exists a block tree $T'$ of $T$ which is an
    \llp-$(1+\e)$-weakly null tree of order \a.
  \item There exists $\gamma<\omega_1$ such that for every \llpkbbt\ $T$
    of order $\gamma$ on a Banach space with a shrinking basis there
    exists a block tree $T'$ of $T$ which is an \llp-$(1+\e)$-block basis
    tree of order \a.
  \end{enumerate}
\end{thm}

This theorem is slightly harder to prove than Theorem~1.1 of~\cite{jo}
because not only do we have to reduce the \llp\ estimate of the nodes, but
we also have to reduce the basis constant to $1+\e$.  In the proof of part
(i), we start off with a weakly null tree of order \b.  We prune the tree
using the Pruning Lemma so that each \sss\ and each branch is a weakly
null basic sequence with basis constant $(1+\e)$.  One may then follow the
same argument as in the proof of~\cite{jo} Theorem~1.1 to reduce an
\llp-$K^2$-weakly null tree isomorphic to $T(\a^2,s)$ to an \llpkwt\ 
isomorphic to $T(\a,s)$, and then complete the proof using the method of
James~\cite{j} as in the last part of the proof of \cite{jo} Theorem~1.1.
Care must be taken to preserve the weakly null structure, but one can
achieve this since $X^*$ is separable and the original trees are weakly
null.

For part (ii) of the theorem we take a block basis tree of order
$\gamma=\ab{\w}{\b}$ which we may assume is isomorphic to the minimal tree
$T(\b,\w)$.  From this we can extract an \llpkbbt\ which is isomorphic to
$T(\b,s)$ with the property that the $s$-subsequences are block bases.
But now the \sss s are weakly null because the basis is shrinking and we
may prune the tree to obtain a tree still isomorphic to $T(\b,s)$ whose
nodes are \llp-$K$ block bases and $(1+\e)$ basic.  We then follow the
same argument as in~\cite{jo} to obtain the result.

\begin{thm}\label{thm:I+wX=w^a}
  If $X$ is a Banach space with separable dual, then $I^+_w(X)=\w^{\a}$
  for some $\a<\w_1$.
\end{thm}
\begin{proof}
  {} From Theorem~\ref{thm:j-reflexive}~(\ref{thm:j-reflexive:item:I}) we
  know that $I^+_w(X)<\w_1$, and so it suffices to show (see e.g.\ 
  Monk~\cite{m}) that if $\b<I^+_w(X)$, then $\ab{\b}{2}<I^+_w(X)$.  We
  may regard $X$ as a subspace of $C[0,1]$, and let $(e_i)_1^{\infty}$ be
  a monotone basis for $C[0,1]$.  Let $T$ be an \llpkwt\ on $X$ of order
  \b.  To make a tree of order \ab{\b}{2} we want to add a tree of order
  \b\ after each terminal node of $T$.

  Let a sequence $(x_i)_1^{\infty}$ in $X$ have property \Q\ for $\e>0$ if
  it is normalized and weakly null.  Let $(x_i)_1^{\infty}$ have property
  \P\ if it is an \e\ perturbation of a normalized block basis of
  $(e_i)_1^{\infty}$.  Properties \P\ and \Q\ clearly satisfy
  conditions PL(1)--(3) of the Pruning Lemma.  Note also that if we apply
  the Pruning Lemma to a tree $T$ for $\e>0$ and a sequence
  $(u_i)_1^l$ satisfying \P, then the resulting sequences
  $(u_1,\dots,u_l,x_1,x_2,\dots)$ are $(1+\e)/(1-\e)$-basic.

  We apply the Pruning Lemma to $T$ with $\e<\min\{1/6,1/(4K)\}$ and the
  empty sequence, so we may assume that for every \sst\
  $(z^i)_1^{\infty}$ of $T$ with $z^i=(x_1,\dots,x_k,y_i)$, the sequence
  $(x_1,\dots,x_k,y_1,y_2,\dots)$ is normalized, weakly null and an $\e$
  perturbation of a normalized block basis of $(e_i)_1^{\infty}$.
  
  Now let $(z^i)_1^{\infty}$ be the sequence of terminal nodes of $T$ with
  $z^i=(x^i_j)_{j=1}^{k_i}$, and apply the Pruning Lemma to $T$ for
  $(x^i_j)_{j=1}^{k_i}$ with $\e<\min\{1/6,1/(4K)\}$ and $\delta=\e$, for
  each $i\geq1$, to obtain a tree $S_i$ which has \P[2\e]\ for the
  sequence $(x^i_j)_{j=1}^{k_i}$.

  To complete the proof we put the trees together as follows.  Let
  \[ S = \bigcup_{i=1}^{\infty}
     \{ (x^i_1),(x^i_1,x^i_2),\dots,
     (x^i_j)_1^{k_i},(x^i_1,\dots,x^i_{k_i},y_1,\dots,y_l) 
     :  (y_j)_1^l\in S_i \} \ .
  \]
  The ordering on $S$ is that inherited from $T$ and the trees $S_i$.  It
  is easy to see that $S$ is an \llpwt.  Indeed, if
  $z=(x_1,\dots,x_k,y_1,\dots,y_l)\in S$, then the sequence is 2-basic and
  both $(x_i)_1^k$ and $(y_j)_1^l$ are \llpk\ sequences so that
  $(x_1,\dots,x_k,y_1,\dots,y_l)$ is an \llp-$6K$ sequence.  Since we have
  extended only the terminal nodes it is clear that the new tree is weakly
  null.  Finally, $o(S)=\ab{\b}{2}$, since $o(S_i)=\b$, and hence
  $(S)^{\b}=T$.
\end{proof}

\begin{thm}\label{thm:I+w<=I+b}
  Let $X$ be a Banach space with shrinking basis $(e_i)_1^{\infty}$.  If
  $I^+_w(X)\geq\w^{\w}$, then $I^+_w(X)=I^+_b(X)$, otherwise, if
  $I^+_w(X)=\w^n$, then $I^+_b(X)=\w^n$ or $\w^{n+1}$.
\end{thm}
\begin{proof}
  We first show by induction on \a\ that if $T$ is an \llpbbt\ on $X$ of
  order \mbox{$\w\!\cdot\!\a$}, then we can extract a certain subtree $S$
  isomorphic to $T(\alpha,s)$.  The tree $S$ will have the property that
  each \sss\ $(x_i)_1^{\infty}$ is a normalized block basis of the
  shrinking basis $(e_i)_1^{\infty}$ and hence is weakly null.  Thus $S$
  will be an \llpwt, and so $I^+_b(X)\leq\ab{\w}{I^+_w(X)}$.  We prove
  this by induction on \a.
  
  For $\a=1$ we may assume that $T\simeq T_{\w}\simeq T(1,\w)$ so that $T$
  consists of disjoint branches $b_n$ of length at least $n$ for each
  $n\geq1$.  Since each branch is a block basis of $(e_i)_1^{\infty}$,
  there must exist a sequence $(n_i)_1^{\infty}$ and a vector $x_i$ in one
  of the nodes of branch $b_{n_i}$ such that $(x_i)_1^{\infty}$ is a
  normalized block basis of $(e_i)_1^{\infty}$.  Set $z^i=(x_i)$ for each
  $i$, then $T'=\{z^i:i\geq1\}$, a sequence of incomparable nodes, is the
  required subtree.
  
  If the result has been proven for \a, then we let $T$ be an \llpbbt\ 
  with order $\ab{\w}{(\a+1)}$ and assume $T\simeq T(\a+1,\w)$, since this
  is a minimal tree of order $\ab{\w}{(\a+1)}$.  Recall that $T(\a+1,\w)$
  is constructed by taking $T_{\w}$, and then after each terminal node
  putting a tree isomorphic to $T(\a,\w)$.  Applying the above argument
  for $\a=1$ to the initial part of the tree, which is isomorphic to
  $T_{\w}$, we may construct a sequence of incomparable nodes
  $(z^i)_1^{\infty}$, with $z^i=(x_i)$ and $(x_i)_1^{\infty}$ a normalized
  block basis of $(e_i)_1^{\infty}$.  After each node $z^i$ we have a tree
  isomorphic to $T(\a,\w)$ from which we may construct an \llp-weakly null
  subtree of order \a\ with the required properties, using the induction
  hypothesis.  Putting these trees together with the nodes
  $(z^i)_1^{\infty}$ we obtain the desired \llpwt.
  
  If \a\ is a limit ordinal and the result has been proven for any ordinal
  smaller than \a, then let $T$ be an \llpbbt\ on $X$ isomorphic to the
  minimal tree $T(\a,\w)$.  Recall that $T(\a,\w)$ is constructed by
  taking a certain sequence of ordinals $(\a_n)_1^{\infty}$ increasing to
  \a\ and letting $T(\a,\w)$ be the disjoint union of the trees $S_n$
  isomorphic to $T(\a_n,\w)$.  By the induction hypothesis we may find a
  subtree $\bar S_n$ of $S_n$ for each $n$ which is isomorphic to
  $T(\a_n,s)$ and has the required properties.  We must now be a little
  careful when putting the trees $\bar S_n$ together.  We cannot just take
  their union since we will need the initial nodes to form a weakly null
  sequence.  Let $(z^{n,i})_{i=1}^{\infty}$ be the sequence of initial
  nodes of $\bar S_n$, with $z^{n,i}=(x_{n,i})$.  Since for each $n$ the
  sequence $(x_{n,i})_{i=1}^{\infty}$ is a normalized block basis of
  $(e_i)_1^{\infty}$, it follows that we may find a sequence
  $(x_i)_1^{\infty}\subset\{x_{n,i} : n,i\geq1\}$ and $1\leq
  k_1<k_2<\cdots$ such that $(x_i)_1^{\infty}$ is a normalized block basis
  of $(e_i)_1^{\infty}$ and $x_i=x_{\phi(i),k_i}$.  Let
  $z^i=(x_i)=z^{\phi(i),k_i}$ and set $S'_n=\{x\in \bar S_n : x\geq z^i,
  \text{ for some } i\in\phi^{-1}(n) \}$.  Then $T'=\cup_nS'_n$ is the
  required tree.  This completes the first part of the proof and shows
  that $I^+_b(X)\leq\ab{\w}{I^+_w(X)}$.
  
  We noted in the proof of
  Theorem~\ref{thm:j-reflexive}~(\ref{thm:j-reflexive:item:I}) that
  $I^+_w(X)\leq I^+_b(X)$, and since we know that both indices are of the
  form $\w^{\a}$ for some $\a<\w_1$ the result follows from the two
  inequalities.
\end{proof}

\section{The Szlenk Index}
\label{sec:szlenk}

In this section we examine the Szlenk index, another isomorphic invariant
of a Banach space, introduced by Szlenk~\cite{sz}.  This is calculated in
a different way to the \llis; it uses collections of subsets in the dual
ball, indexed by countable ordinals.  We show that the Szlenk index is in
fact the same as the \llpwi\ provided the space does not contain $\ell_1$.

\begin{defn}\label{defn:szlenk}
  For a fixed $\e>0$ we construct inductively sets $P_{\a}(\e)\subseteq
  B_{X^*}$.  Let $P_0(\e)=B_{X^*}$ and if we have constructed
  $P_{\a}(\e)$, then let
  \[  
    P_{\a+1}(\e) = \{ f\in B_{X^*} : \exists(f_m)_1^{\infty} \subset
    P_{\a}(\e) \mbox{ with } f_m\stackrel{w^*}{\rightarrow}f \mbox{ and } 
    \lim\inf\|f_m-f\|\geq\e \} \ . 
  \]  
  If \a\ is a limit ordinal and we have chosen $P_{\b}(\e)$ for each
  $\b<\a$, then let
  \[ P_{\a}(\e) = \bigcap_{\b<\a}P_{\b}(\e) \ . \]
  We define the \emph{\e-Szlenk index} of $X$ to be 
  \[ \eta(\e,X) = \sup\{ \a : P_{\a}(\e)\neq\emptyset \} \ , \]
  if such an $\alpha$ exists, and $\omega_1$ otherwise, and the
  \emph{Szlenk index} of a Banach space $X$ as
  \[ \eta(X) = \sup_{\e>0}\eta(\e,X) \ . \]
  One can show that if \llo\ does not embed into $X$, then
  $\eta(X)<\omega_1$ if and only if $X^*$ is separable.  Indeed, if $X^*$
  is not separable, then by Stegall's result used above~\cite{st} there
  exists a homeomorphic copy of $\Delta$ in $(B_{X^*},w^*)$ which is
  $(1-\e)$-separated, i.e.\ $\|x^*-y^*\|>1-\e$ for $x^*, y^*\in\Delta$
  with $x^*\neq y^*$.  Thus $\Delta\subseteq P_{\alpha}(1/2)$ for each
  $\alpha<\omega_1$.
\end{defn}

In his original definition Szlenk used sets $P'_{\alpha}(\epsilon)$
defined in a similar manner to the sets $P_{\a}(\e)$ above, except at
successor ordinals he had
\begin{multline*}
  P'_{\a+1}(\e) = \{ f\in X^* : \exists(x_m)_1^{\infty}\subset B_X,
  (f_m)\subset P'_{\a}(\e) \\
  \mbox{such that } f_m\stackrel{w^*}{\rightarrow}f,
  x_m\stackrel{w}{\rightarrow}0 \mbox{ and } \lim\sup|f_m(x_m)|\geq\e\}\ .
\end{multline*}
The original Szlenk indices $\eta'(\e,X)$ and $\eta'(X)$ were defined as
before, but using the sets $P'_{\alpha}(\epsilon)$.  Szlenk then showed
that if $X^*$ is separable, then $\eta'(X)<\omega_1$.  These two
definitions may give different values for the \e-indices.  However, we
have using Rosenthal's $\ell_1$ theorem that if $X$ does not contain
$\ell_1$, then
\[ P'_{\a}(\e) \subseteq P_{\a}(\e) \subseteq P'_{\a}(\e/2) \]
and hence 
\[ \eta'(\e,X) \leq \eta(\e,X) \leq \eta'(\e/2,X) \ . \]
Thus, if $\ell_1\not\hookrightarrow X$, then $\eta'(X)=\eta(X)$.  Since we
shall only be considering spaces which do not contain $\ell_1$, in the
sequel we shall apply the definition for the Szlenk index using the sets
$P_{\a}(\e)$.

\begin{thm}\label{thm:sz=i+w}
  If $X$ is a separable Banach space not containing \llo, then
  $\eta(X)=I^+_w(X)$. 
\end{thm}
This result shows that despite the Szlenk index being calculated using
subsets of the unit ball of $X^*$, while the \llpwi\ is calculated using
trees on the unit ball of $X$, the two indices are in fact the same.
Moreover we show that from the sets $P_{\a}(\e)$ used to calculate the
Szlenk index one can generate trees which are analogous to the trees used
in the \llpwi.  

By Theorem~\ref{thm:j-reflexive}~(\ref{thm:j-reflexive:item:I}) and the
remarks on the Szlenk index above we have that if \llo\ does not embed
inside $X$, $\eta(X)<\omega_1$ if and only if $X^*$ is separable if and
only if $I^+_w(X)<\omega_1$.  Thus in proving the theorem we may restrict
ourselves to the case where $X^*$ is separable.  The proof is in two
parts.  In the first we show that if $P_{\a}(\e)\neq\emptyset$, then there
exists an \llpwt\ on $X$ of order \a\ with constant $8/\e$.  In the second
part we demonstrate that if we have an \llpwt\ on $X$ with constant $K$
and order $\a$, then $P_{\a}(1/K)\neq\emptyset$.  Thus our first task is
to prove

\begin{prop}\label{prop:sz=>i+w}
  If $P_{\a}(\e)\neq\emptyset$, then there exists an \llpwt\ on $X$ of
  order \a\ with constant $8/\e$.
\end{prop}
To prove this proposition we first construct a tree of order \a\ 
isomorphic to $T(\a,s)$ on $B_{X^*}$.  From this tree we construct an
isomorphic tree on $B_X$ which is an \llpwt.  We construct the tree on
$B_{X^*}$ in the next two lemmas, and in Lemma~\ref{lem:C} describe the
properties \P\ and \Q\ needed to construct the tree on $B_X$ from the tree
on $B_{X^*}$.

\begin{lem}\label{lem:A}
  If $P_{\a}(\e)\neq\emptyset$, then for each $f_0\in P_{\a}(\e)$ and each
  \wstar\ relatively open neighborhood $O$ of $f_0$, with respect to
  $P_{\a}(\e)$, there exists a tree $T$ on $O$, isomorphic to $T(\a,s)$,
  such that if $z=(f_i)_1^k\in T$ has immediate successors
  $(z^j)_1^{\infty}$ with $z^j=(f_1,\dots,f_k,g_j)$, then $g_j\wkstar f_k$
  as $(j\rightarrow\infty)$ and $\lim\inf_j\|g_j-f_k\|\geq\e$.  Also, if
  $(z^j)_1^{\infty}$ is the sequence of initial nodes, with $z^j=(g_j)$,
  then $g_j\wkstar f_0$ $(j\rightarrow\infty)$ and
  $\lim\inf_j\|g_j-f_0\|\geq\e$.
\end{lem}

\begin{proof}[Proof of Lemma~\ref{lem:A}]
  As usual we use induction on \a.  For the initial case $\a=1$, let
  $f_0\in P_1(\e)$ and let $O$ be a \wstar\ relatively open neighborhood
  of $f_0$.  We may find $(g_j)_1^{\infty}\subset O$ such that $g_j\wkstar
  f_0$ and $\lim\inf\|g_j-f_0\|\geq\e$.  Set $z^j=(g_j)$, then
  $T=\{z^j:j\geq1\}$ is the required tree.
  
  We next suppose the result has been proven for \a; let $f_0\in
  P_{\a+1}(\e)$, and let $O$ be a \wstar\ relatively open neighborhood of
  $f_0$.  We may find $(g_j)_1^{\infty}\subset O$ with $g_j\wkstar f_0$
  and $\lim\inf\|g_j-f_0\|\geq\e$.  From the induction hypothesis, for
  each $j\geq1$ there exists a tree $S_j$ isomorphic to $T(\a,s)$
  satisfying the requirements for $g_j\in P_{\a}(\e)$.  Let $S'_j =
  \{(g_j,h_1,\dots,h_k) : (h_i)_1^k\in S_j\}$ so that the trees $S'_j$ are
  disjoint.  Define $T=\cup_jS'_j$, then $T\simeq T(\alpha+1,s)$ and
  satisfies the requirements of the lemma for $f_0$.
  
  If \a\ is a limit ordinal and the result has been proven for each
  $\b<\a$, let $(\a_n)_1^{\infty}$ be the sequence of successor ordinals
  increasing to \a\ so that $T(\a,s)=\cup_nT(\a_n,s)$.  Let $f_0\in
  P_{\a}(\e)$ (so that $f_0\in P_{\a_n}(\e)$ for each $n\geq1$) and let
  $O\supseteq O_1\supseteq O_2\supseteq\cdots$ be a decreasing collection
  of \wstar\ relatively open neighborhoods of $f_0$, so that
  $\cap_iO_i=\{f_0\}$, which may be chosen since $X$ is separable.  By the
  induction hypothesis we may find a tree $S_n$ isomorphic to $T(\a_n,s)$
  satisfying the lemma for $f_0$ and $O_n$ for each $n\geq1$.  Let
  $(z^{n,i})_{i=1}^{\infty}$ be the sequence of initial nodes of $S_n$
  with $z^{n,i}=(f_{n,i})$, so that $f_{n,i}\wkstar f_0\ 
  (i\rightarrow\infty)$ and $f_{n,i}\in O_n$ for each $i\geq1$ and every
  $n\geq1$.  Since the sets $O_n$ are decreasing we may find a subsequence
  $(f_i)_1^{\infty}$ of $\{f_{n,i}:n,i\geq1\}$ and numbers $1\leq
  k_1<k_2<\cdots$ with $f_i=f_{\phi(i),k_i}$ for each $i$ (where $\phi$ is
  the function from the Pruning Lemma) such that $f_i\wkstar f_0$ and
  $\lim\inf\|f_i-f_0\|\geq\e$.  It is at this limit ordinal stage that we
  use the relatively open neighborhoods $O_n$ to ensure that we choose
  $f_i\in S_i$, so that the order of the tree $T$ below will be $\alpha$.
  The tree
  \[ T=\bigcup_n \{ z\in S_n : z\geq f_{i,k_i}, i\in\phi^{-1}(n) \} \]
  satisfies the requirements of the lemma.
\end{proof}

\begin{lem}\label{lem:B}
  If $P_{\a}(\e)\neq\emptyset$, then for any $\delta>0$ there exists a
  tree $T$ isomorphic to $T(\a,s)$ on $B_{X^*}$ such that
  \begin{enumerate}[\rm (i)]
  \item $h_j\wkstar0$ and $\e/2-\delta\leq\|h_j\|\leq1\ (j\geq1)$ for
    every \sss\ $(h_j)_1^{\infty}$ of $T$;
  \item $\|\sum_k^lg_i\|\leq1\ (1\leq k\leq l\leq m)$ for every
    $(g_i)_1^m\in T$.
  \end{enumerate}
\end{lem}
\begin{proof}
  Let $f_0\in P_{\a}(\e)$ and let $T$ be the tree for $f_0$ from the
  previous lemma for $O=P_{\alpha}(\varepsilon)$.  Replace each node
  $z=(f_i)_1^m\in T$ with the node
  \[  
  \bar z = ({\textstyle \frac{1}{2}}(f_1-f_0), 
  {\textstyle \frac{1}{2}}(f_2-f_1),\dots,
  {\textstyle \frac{1}{2}}(f_m-f_{m-1}))  
  \]
  to obtain the tree $\bar T$ which is still isomorphic to $T(\a,s)$.
  Clearly, if $(h_j)_1^{\infty}$ is any \sss, then $h_j\wkstar0$ and
  $\lim\inf_j\|h_j\|\geq\e/2$.  Let $(x^*_i)_1^{\infty}$ have property \Q\ 
  if it is weak$^*$ null with $\lim\inf_j\|x^*_j\|\geq\e/2$, and property
  \P\ if it is weak$^*$ null with $\|x^*_j\|\geq\e/2$ for every $i$.  It is
  clear that \Q\ and \P\ satisfy condition PL(1) of the Pruning Lemma, as
  modified by Remark~\ref{rem:pruning}~(iii) after it and we obtain
  condition PL(2) using the fact that $X$ is separable.  Thus we may apply
  the Pruning Lemma and prune $\bar T$ to obtain a tree $T'$ with property
  \P[\e-\delta] satisfying condition (i) above.  To see that (ii) holds,
  note that each node $z=(g_i)_1^m\in T$ is of the form
  $(\frac{1}{2}(f_1-f_0), \frac{1}{2}(f_2-f_1), \dots,
  \frac{1}{2}(f_m-f_{m-1}))$ so that
  $\|\sum_k^lg_i\|=\frac{1}{2}\|f_l-f_{k-1}\|\leq1$.
\end{proof}

Our next lemma contains the basic relationship between the \wstar\ null
trees constructed above and trees in $X$.  The lemma is stated so as to
verify the hypotheses of the Pruning Lemma with the additional conditions
of Remark~\ref{rem:pruning}~(iv).  

\begin{lem}\label{lem:C}
  Let $X$ be a Banach space with separable dual, let $\e>0$ and let
  $(f_i)_1^{\infty}\subseteq B_{X^*}$ be weak$^*$ null with
  $\e/2\leq\|f_i\|\leq1$ for every $i\geq1$.  Then there exists a
  subsequence $(f'_i)_1^{\infty}$ of $(f_i)_1^{\infty}$ and a weakly null
  sequence $(x_i)_1^{\infty}\subseteq S_X$ such that $f'_i(x_i)\geq\e/5$,
  $|f'_i(x_j)|<\e/2^{i+6}$ whenever $i\neq j$ and $(x_i)_1^k$ is
  $(1+\e(1-2^{-k}))$ basic for every $k\geq1$.
\end{lem}

\begin{proof}
  Let $0<\delta<\epsilon$, to be chosen later.  We first choose a sequence
  $(y_i)_1^{\infty}\subset S_X$ with $f_i(y_i)\geq\e/2-\delta$ for each
  $i$.  Since $X^*$ is separable we may assume $(y_i)_1^{\infty}$ is
  weakly Cauchy (by taking a subsequence of $(y_i)_1^{\infty}$ and then
  the same subsequence of $(f_i)_1^{\infty}$).  Again, by taking
  subsequences and using that $f_i\wkstar0$, we may assume that
  $|f_n(y_i)|<\theta(\e,n)$ if $i<n$ (where $\theta(\e,n)$ is small, to be
  chosen later).  Now set
  \[ x_n=\frac{y_n-y_{n-1}}{\|y_n-y_{n-1}\|} \ , \]
  so that, since $f_n(y_n-y_{n-1})\geq\e/2-\delta-\theta(\e,n)$, it
  follows that $\|y_n-y_{n-1}\|\geq\e/2-\delta-\theta(\e,n)$, and hence
  $x_n\wkly0$ and $f_n(x_n)\geq\e/4-\delta/2-\theta(\e,n)/2$.  Further,
  for $i<n$,
  \[
    |f_n(x_i)| = \frac{|f_n(y_i)-f_n(y_{i-1})|}{\|y_i-y_{i-1}\|}
               < \frac{2\theta(\e,n)}{\e/2-\delta-\theta(\e,n)} \ .
  \]
  If $\delta<\e/40$ and $\theta(\e,n)=\e^2/2^{i+10}$, then
  $f_i(x_i)\geq\e/5$ and $|f_n(x_i)|<\e/2^{n+6}$ when $i<n$.  Next, since
  $x_n\wkly0$, we may pass to subsequences of $(f_i)_1^{\infty}$ and
  $(x_i)_1^{\infty}$ to obtain $|f_i(x_j)|<\e/2^{i+6}$ when $i\neq j$.  We
  now pass to one last pair of subsequences $(f'_i)_1^{\infty}$ and
  $(x'_i)_1^{\infty}$ so that $(x'_i)_1^k$ is $(1+\e(1-2^{-k}))$ basic
  for every $k\geq1$ as required.
\end{proof}

\begin{proof}[Proof of Proposition~\ref{prop:sz=>i+w}]
  We have that $P_{\a}(\e)\neq\emptyset$ and we want to construct an
  \llpwt\ on $X$ of order \a\ and constant $K=K(\e)=8/\e$.  Let $T$ be the
  tree on $B_{X^*}$ for some $\delta<\e/12$ from Lemma~\ref{lem:B}.  We
  want to construct a tree $S$ in $S_X$, isomorphic to $T$, so that if
  $(f_i)_1^m\in T$ has immediate successors $(f_1,\dots,f_m,f_{m+1}),
  (f_1,\dots,f_m,f_{m+2}),\dots$ etc., and $\xim, (x_1,\dots,x_m,x_{m+1}),
  (x_1,\dots,x_m,x_{m+2}),\dots$ are the corresponding nodes of $S$, then
  $(x_i)_1^{\infty}$ is weakly null, $f_i(x_i)\geq\e/5$,
  $|f_i(x_j)|<\e/2^{i+6}$ whenever $i\neq j$ and $(x_i)_1^k$ is
  $(1+\e(1-2^{-k}))$ basic for every $k\geq1$.  The proof is very similar
  to that of the Pruning Lemma, although a little stronger as we must keep
  track of two trees $S$ and $T$, so we will not give it here.
  
  We claim that $S$ is the required \llpkwt.  We already know that $S$ is
  a weakly null tree.  We must show that if $(x_i)_1^m\in S$, then
  $(x_i)_1^m$ is an \llpk-sequence.  We know that $(x_i)_1^m$ is $(1+\e)$
  basic, so we seek $f\in S_{X^*}$ such that $f(x_i)\geq8/\e$ for $1\leq
  i\leq m$.  Let $(f_i)_1^m$ be the corresponding node in $T$ to \xim\ and
  recall that $\|\sum_1^mf_i\|\leq1$.  Now,
  \[
    \sum_{j=1}^m f_j(x_i)  = f_i(x_i) + \sum_{j\neq i}f_j(x_i) 
                           \geq f_i(x_i) - \sum_{j\neq i}|f_j(x_i)| 
                           \geq \frac{\e}{6} - \sum_{j\neq i}2^{-j-6}\e 
                           \geq \frac{\e}{8} \ ,
  \]
  for $1\leq i\leq m$.  Finally, setting
  $f=\sum_{j=1}^mf_j/\|\sum_{j=1}^mf_j\|$ we still have $f(x_i)\geq\e/8$
  for each $i$, and hence $(x_i)_1^m$ is an \llpk-sequence with $K=8/\e$.
\end{proof}

\begin{rem}\label{rem:sz=>i+w}
  One can be more careful with the estimates in the proofs of
  Lemma~\ref{lem:C} and Proposition~\ref{prop:sz=>i+w}, and obtain  an
  \llpkwt\ on $X$ of order \a\ with $1/K=\e/4-\delta$ for any $\delta>0$.
\end{rem}

This completes the first part of the proof of Theorem~\ref{thm:sz=i+w}.
We now have to show how to get from an \llpwt\ on $B_X$ to the sets
required in the calculation of the Szlenk index.

\begin{defn}\label{defn:branch-fnl}
  If $T$ is an \llpkt\ on $X$ and $(x_i)_1^n$ is a terminal node of $T$,
  then let $\gamma=\{ y\in T : y\leq(x_i)_1^n \}$ be the branch of $T$
  ending at $(x_i)_1^n$.  A $K$-\emph{branch functional of $\gamma$} is an
  element $f_{\gamma}\in KB_{X^*}$ with $f_{\gamma}(x_i)\geq1$ for each
  $i$.  These exist from the equivalent formulation of \llp\ sequences in
  Fact~\ref{rem:l1p-equiv}.  A \emph{full set of branch functionals of
    $T$} is a subset of $X^*$ which contains a branch functional for each
  branch of $T$.
\end{defn}

\begin{lem} 
  If $T$ is an \llpkwt\ of order \a\ and $W$ is a weak$^*$ closed subset of
  $KB_{X^*}$ which contains a full set of $K$-branch functionals of $T$,
  then $W\cap KP_{\a}(1/K)\neq\emptyset$.
\end{lem}
\begin{proof}
  As usual we proceed by induction on \a.  If $o(T)=1$, then
  $T=\{z^i:i\geq1\}$ where $z^i=(x_i)$ and $(x_i)_1^{\infty}$ is a
  normalized weakly null sequence.  For each $i$ pick $f_i\in W$ with
  $f_i(x_i)\geq1$, then choose a subsequence $(f_{n_i})_1^{\infty}$ which
  converges weak$^*$ to some $f\in W$.  Choose a sequence $\e_i\searrow0$;
  since $(x_i)_1^{\infty}$ is weakly null, it follows that for each $i$
  there exists $m_i\geq1$ such that $|f(x_j)|<\e_j$ for every $j\geq m_i$.
  But now
  \[ 
  |f_{n_j}(x_{n_j})-f(x_{n_j})| \geq1-\e_i \text{ for every } j\geq m_i \
   , 
  \] 
  so that $\lim\inf\|f_{n_j}-f\|\geq1$ and hence $f\in KP_1(1/K)$.
  The result for the case $\a=1$ follows easily from this.
  
  If the result has been proven for \a, let $T$ be an \llpkwt\ of order
  $\a+1$ and let $W$ be a weak$^*$ closed subset of $KB_{X^*}$ which
  contains a full set of branch functionals of $T$.  Let
  $(z^i)_1^{\infty}$ be the sequence of initial nodes of $T$ with
  $z^i=(x_i)$ and $(x_i)_1^{\infty}$ a normalized weakly null sequence.
  For each $i$ let $W_i$ be the weak$^*$ closure of the set of branch
  functionals of $T$ in $W$ for branches whose initial node is $z^i$.
  Thus $f(x_i)\geq1$ for every $f\in W_i$.  Further, let $T_i=\{ y\in T:
  y>z^i\}$, so that $T_i$ is an \llpkwt\ of order \a, and $W_i$ is a
  weak$^*$ closed subset of $KB_{X^*}$ which contains a full set of branch
  functionals of $T_i$.  Hence, by the induction hypothesis, there exists
  $f_i\in W_i$ with $f_i\in W_i \cap KP_{\a}(1/K)$.  Now, $f_i(x_i)\geq1$
  for each $i$, so we may now proceed as in the case $\a=1$ to obtain
  $f_{n_i}\wkstar f$ with $\lim\inf\|f_{n_i}-f\|\geq1$.  Thus $f\in
  KP_{\a+1}(1/K)$. Then, since $W$ is weak$^*$ closed, and since
  $W_i\subseteq W$ for each $i$, it follows that $f\in W\cap
  KP_{\a+1}(1/K)$ as required.
  
  For the case where \a\ is a limit ordinal we simply note that if the
  result has been proven for each $\b<\a$, and if we have $T$ and $W$ as
  in the statement of the lemma, then $ W\cap KP_{\b}(1/K)\neq\emptyset$
  for each $\b<\a$.  This forms a countable decreasing sequence of
  non-empty weak$^*$ closed sets in the weak$^*$ compact set $KB_{X^*}$.
  Thus $ W\cap KP_{\a}(1/K)= W\cap(\cap_{\b<\a}KP_{\b}(1/K))=\cap_{\b<\a}(
  W\cap KP_{\b}(1/K))\neq\emptyset$, which completes the proof.
\end{proof}

\begin{prop}\label{prop:sz<=i+w}
  If there exists an \llpkwt\ on $X$ of order \a, then
  $P_{\a}(1/K)\neq\emptyset$.
\end{prop}
\begin{proof}
  If $T$ is an \llpkwt\ on $X$ of order \a, then there exists a branch
  functional for each branch of $T$.  Thus we may take $W=KB_{X^*}$ in the
  above lemma, to obtain $P_{\a}(1/K)= \frac{1}{K}W\cap
  P_{\a}(1/K)\neq\emptyset$.
\end{proof}

\section{The $\ell_1$ index of the Schreier spaces and the $C(\alpha)$
  spaces} 
\label{sec:Xa-Ca-indices}

In this section we calculate the \llis\ of the Schreier spaces and the
$C(\alpha)$ spaces using the results from the previous two sections.  We
first give some notation.

\begin{defn}[\cite{aa}]\label{defn:schreier-sets}
  Let $E,F$ be subsets of $\N$ and $n\geq1$.  We write $E<F$ if $F$ is
  empty or $\max E<\min F$; we write $n<E$ if $\{n\}<E$, and $n\leq E$ if
  $n=\min E$ or $n<E$.  The Schreier sets $\cs$, for each
  $\alpha<\omega_1$, are defined inductively as follows: Let
  $\cs[0]=\{\{n\} : n\geq1\}\cup\{\emptyset\}$ and $\cs[1]=\{F\subset
  {\mathbf N} : |F|\leq F\}$.  (Note that this definition allows for
  $\emptyset\in\cs[1]$.)  If $\cs$ has been defined, let
  \[\cs[\alpha+1]=\{ \cup_1^kF_i : k\leq
  F_1<\dots<F_k,\ F_i\in \cs\ (i=1,\dots,k),\ k\in\N \}\ .\] If
  $\alpha$ is a limit ordinal with $\cs[\beta]$ defined for each
  $\beta<\alpha$, choose and fix an increasing sequence of ordinals
  $(\alpha_n)$ with $\alpha=\sup_n\alpha_n$ and let
  \[\cs=\bigcup_{n=1}^{\infty}\{ F\in\cs[\alpha_n] : n\leq F\}\
  .\] 
\end{defn}
  
Each $\cs$ has the following two important properties.  First, if
$F=\{m_1,\dots,m_k\}\in\cs$ and $n_1<\dots<n_k$ satisfies: $m_i\leq
n_i$ for $i\leq k$, then $\{n_1,\dots,n_k\}\in\cs$ (this is called
\emph{spreading}).  Second, whenever $E\subset F$ and $F\in\cs$ then
$E\in\cs$ (this is called \emph{hereditary}).

For each $\alpha<\omega_1$ the Schreier set $\cs$ generates a
tree, ${\rm Tree}(\cs)=(\cs,\subseteq)$, ordered by
inclusion.  It is easy to see that the order of ${\rm Tree}(\cs)$ is
$\omega^{\alpha}+1$~\cite{aa}.

\begin{defn}\label{defn:x-alpha}
  The Schreier spaces generalize Schreier's example~\cite{sh}; they were
  introduced in~\cite{ao} for \a\ finite and in~\cite{aa} for \a\ 
  infinite.  We first define $c_{00}$ to be the linear space of all real
  sequences with finite support, and let $(e_i)_1^{\infty}$ be the unit
  vector basis of $c_{00}$.  For each $\alpha<\omega_1$ let
  \mbox{$\|\cdot\|_{\alpha}$} be the norm on $c_{00}$ given by:
  \[ \left\|\sum a_ie_i\right\|_{\alpha}= \sup_{E\in{\mathcal 
      S}_{\alpha}}\Bigl|\sum_{i\in E}a_i\Bigr| \ , \] 
  then \emph{the Schreier space $X_{\alpha}$} is the completion of
  $(c_{00},\|\cdot\|_{\alpha})$.  Note that because $\cs$ is
  hereditary $(e_i)_1^{\infty}$ is a normalized 1-unconditional basis for
  $X_{\alpha}$.  
\end{defn}

\begin{defn}\label{defn:c-alpha}
  If $\alpha>0$ is an ordinal, then $C(\alpha)$ denotes the Banach space
  of all continuous real-valued functions on the ordinals less than or
  equal to $\a$, where $[1,\a]=\{\b : \b\leq\a\}$ has the order topology,
  with the norm $\|x\|=\sup_{\beta\in\alpha+1}|x(\beta)|$.  Thus
  $C(\alpha)=C([1,\a])$ is the space of all continuous functions
  $x:[1,\a]\rightarrow\R$.
\end{defn}

The following classical theorem of Bessaga and Pe\l czy\'nski~\cite{bp}
partitions the $C(\alpha)$ spaces $(\omega\leq\alpha<\omega_1)$ into
isomorphism classes.
\begin{thm}[Bessaga and Pe\l czy\'nski] \label{thm:bp}
  Let $\w\leq\a\leq\b<\w_1$, then $C(\alpha)$ is isomorphic to $C(\beta)$
  if, and only if, $\b<\a^{\w}$.  Furthermore, if we do have $\b<\a^{\w}$,
  then $C(\beta)\oplus C(\alpha)$ is isomorphic to $C(\alpha)$.
\end{thm}
Thus, in studying isomorphic invariants of the spaces $C(\alpha)$ for
$\alpha<\omega_1$, and hence in particular when calculating the \llis, it
suffices to consider the spaces $C(\omega^{\omega^{\beta}})$ for
$\beta<\omega_1$.  It is well known (see \cite{ab}) and not difficult to
see that the Szlenk indices of the $C(\a)$ spaces are given by
$\eta(\C)=\w^{\a+1}$ and so by Theorem~\ref{thm:sz=i+w}
$I^+_w(\C)=\w^{\a+1}$.

The main result of this section is the following theorem.

\begin{thm}\label{thm:Xa-Ca-indices}
  For $1\leq\alpha<\omega_1$
  \begin{enumerate}
  \item $I(X_{\alpha})=I_b(X_{\alpha})=\omega^{\alpha+1}$ with respect to
    the unit vector basis \seq{e_i} of $X_{\a}$;
  \item $I_b(\C)=\omega^{\alpha+1}$ with respect to the node basis,
    described below;
  \item $I(\C)=\omega^{1+\alpha+1}$.
  \end{enumerate}
\end{thm}
Notice that if $\a=n$ is finite, then $I_b(\C[n])<I(\C[n])$.  Note also
that since the unit vector basis for $X_{\a}$ is unconditional, it follows
that $I^+_b(X_{\a})=I_b(X_{\a})$.  Neither $X_{\a}$, nor \C\ is reflexive,
so $I^+(X_{\a})=I^+(\C)=\w_1$.

In order to prove the above theorem we must first describe the node basis
for \C\ and clarify the relationship between the Schreier sets, the
Schreier spaces $X_{\a}$ and \C.

We know that if we identify \cs\ with $\{\mathbf{1}_F :
F\in\cs\}\subset\{0,1\}^{\N}$, then \cs\ is homeomorphic to
$[1,\w^{\w^{\a}}]$ in the topology of pointwise convergence (see \cite{aa}
and \cite{ms}).  Thus we shall consider this representation of
$[1,\w^{\w^{\a}}]$ in the sequel.  

Define a partial order on \cs\ by $F\prq G$ if and only if $F$ is an
initial segment of $G$, i.e.\ $G\cap\{1,2,\dots\max F\}=F$.  This order
induces a natural tree structure on \cs.  For each $F\in\cs$ define a
function $\chir_F:\cs\rightarrow\{0,1\}$ by
\begin{equation*}
  \chir_F(G)=
  \begin{cases}
    1, & \text{if $F\prq G$;}  \\
    0, & \text{otherwise.}
  \end{cases}
\end{equation*}
This function is thus 1 on every $G$ in \cs\ that extends $F$.  Let
$\B=\{\chir_F:F\in\cs\}$, then we say that $(\chir_{F_i})_{i=0}^{\infty}$
is an admissible enumeration of \B\ if and only if $F_i\prec F_j$ implies
$i<j$.  Since $\emptyset\prq F$ for every $F\in\cs$, it follows that
$F_0=\emptyset$ in any admissible enumeration of \B.  Notice that
admissible enumerations preserve the tree structure in that we have an
order preserving map from $(\cs,\prec)$ to \N.  Furthermore, for each
$F\in\cs$ we have $\chir_F\in\C$.  Indeed, if $(G_i)_{i=1}^{\infty}$ is a
sequence in \cs\ converging to $G$, then for every $N\geq1$ there exists
$I\geq1$ such that $G_i\cap\{1,\dots,N\}=G\cap\{1,\dots,N\}$ for every
$i\geq I$, and in particular there exists $I_G\geq1$ such that
$G_i\cap\{1,\dots,\max G\}=G$ for every $i\geq I_G$.  It is now clear that
$F\prq G$ if and only if $F\prq G_i$ for every $i\geq I_G$ and hence
$\chir_F\in\C$ as required.

\begin{lem}\label{lem:mono-basis}
  If $(\chir_{F_i})_{i=0}^{\infty}$ is an admissible enumeration of \B,
  then $(\chir_{F_i})_{i=0}^{\infty}$ is a monotone basis for \C.
\end{lem}
\begin{proof}
  We first show that $(\chir_{F_i})_{i=0}^{\infty}$ is a monotone basic
  sequence, and then apply the Stone-Weierstrass theorem to obtain that
  its span is all of \C.

  Let $(a_i)_0^{\infty}\in\R$, then 
  \[
    \Bigl\| \sum_{i=0}^{k+1} a_i\chir_{F_i} \Bigr\| 
      = \sup_{G\in\cs} \Bigl| \sum_{i=0}^{k+1} a_i\chir_{F_i}(G) \Bigr| 
      = \sup_{G\in\cs} \Bigl| \sum_{\substack{i:F_i\prq G \\ 0\leq i\leq
                                     k+1}}a_i \Bigr|\ . 
  \]
  Since \cs\ is hereditary and the enumeration of \B\ is admissible, it
  follows that there is an index $i_0\leq k$ and $j>F_{i_0}$ such that
  $F_{k+1}=F_{i_0}\cup\{j\}$.  Next, observe that for all $G\in\cs$ with
  $F_{k+1}\prq G$, if $G'=G\setminus\{j\}$, then for $i\leq k+1$, $F_i\prq
  G'$ if and only if $F_i\prq G$ and $i\leq k$ (since $G\cap F_{k+1} =
  G\cap\{1,\dots,j\} = F_{k+1}$).  Thus
  \[ \sup_{G\in\cs} \Bigl| \sum_{\substack{i:F_i\prq G \\ 0\leq i\leq
                                     k+1}} a_i \Bigr| \geq
     \sup_{G\in\cs} \Bigl| \sum_{\substack{i:F_i\prq G \\ 0\leq i\leq
                                     k}} a_i \Bigr| =
     \Bigl\| \sum_{i=0}^k a_i\chir_{F_i} \Bigr\| \ , \]
  and so $(\chir_{F_i})_{i=0}^{\infty}$ is a monotone basic sequence.
  
  To see that $[\chir_F:F\in\cs]=\C$ we shall apply the Stone-Weierstrass
  theorem.  Since $\chir_{\emptyset}=1$ for each $G\in\cs$, it follows
  that $[\chir_F:F\in\cs]$ contains the constant function.  It is
  easy to see that $[\chir_F:F\in\cs]$ separates the points of \cs, so it
  remains to show that the set contains the algebra generated by
  $\{\chir_F:F\in\cs\}$.  If $F,F'\in\cs$ and $\chir_F\cdot\chir_{F'}$ is
  not identically zero, then there exists $G\in\cs$ such that
  $\chir_F(G)=\chir_{F'}(G)=1$, i.e., $F\prq G$ and $F'\prq G$ so that
  both $\{1,2,\dots,\max F\}\cap G = F$, and $\{1,2,\dots,\max F'\}\cap G
  = F'$.  Hence either $F'\prq F$, or $F\prq F'$ which gives
  $\chir_F\cdot\chir_{F'}$ is $\chir_F$ or $\chir_{F'}$ respectively.  In
  either case we have that the algebra is contained in the linear span, as
  required, which completes the proof.
\end{proof}

\begin{defn}\label{defn:adm-enum}
  Since any admissible enumeration of \B, is a monotone basis for \C, we
  shall call \B\ the \emph{node basis} for \C.
\end{defn}

\begin{rem}\label{rem:shrinking}
  For any point $\b\in[1,\w^{\w^{\a}}]$ there are only finitely many
  elements in the node basis which have \b\ in their support.  With this
  in mind it is clear that the node basis is shrinking.
\end{rem}

Finally let us consider the spaces $X_{\a}$.  We have defined 
\[ \Bigl\|\sum a_ie_i\Bigr\|_{\alpha}= \sup_{G\in{\mathcal
    S}_{\alpha}}\Bigl|\sum_{i\in G}a_i\Bigr| \ . \] For each $i\geq1$ let
$f_i=\ch_{\{G\in\cs:i\in G\}}$.  Clearly $f_i\in\C$ for each $i$, and in
\C:
\[
  \Bigl\|\sum a_if_i\Bigr\|_{\C} 
         = \sup_{F\in\cs} \Bigl| \sum a_if_i(F) \Bigr| 
         = \sup_{F\in\cs} \Bigl| \sum_{i\in F} a_i \Bigr| 
         = \Bigl\| \sum a_ie_i \Bigr\|_{\alpha} \ . 
\]
Thus \seq{f_i} is 1-equivalent to \seq{e_i} and $X_{\a}$ can be
isometrically embedded in \C.  Actually, more is true.

\begin{lem}\label{lem:Xa-bb-Cwwa}
  If $(\chir_{F_i})_0^{\infty}$ is an admissible enumeration of \B, then
  the basis $(f_i)_1^{\infty}$ 
  for $X_{\a}$ is 1-equivalent to a block basis of
  $(\chir_{F_i})_0^{\infty}$.
\end{lem}
\begin{proof}
  The heart of the proof lies in choosing an appropriate block basis of
  $(\chir_{F_i})_0^{\infty}$.  To do this we shall construct a tree
  isomorphism $\psi$ from $\cs$ into $\cs$ by induction, and then the
  map from the basis $(f_i)_1^{\infty}$ of $X_{\a}$ to a block basis of
  $(\chir_{F_i})_0^{\infty}$ will be given by 
  \[ U(f_i) = \sum_{\substack{G\in\cs \\ \max G=i}} \chir_{\psif G} \
  . \] 
  This immediately gives the ordering requirement on $\psi$ that if $\max
  G=n$ and $\max G'=n+1$, then $\chir_{\psif G}$ precedes
  $\chir_{\psif{G'}}$, i.e.\ if $\psi G=F_i$ and $\psi G'=F_j$, then
  $i<j$. 
  
  To help us write down the construction of $\psi$ more explicitly we
  define subtrees $T_F$ of \cs, for $F\in\cs$, by $T_F = \spp \chir_F = \{
  G\in\cs : F\prq G\}$, with the order $\prq$, and as usual $o(T_F)$ is
  the order of the tree.  Clearly $T_{\emptyset}=\cs$ and $o(T_F)$ is a
  successor ordinal for each $F\in\cs$, since $F$ is the unique initial
  node.
  
  Let $F$ be a non-terminal node of the tree $T_{\emptyset}$, so that
  there exists $G\in\cs$ with $F\prec G$.  Then there are infinitely many
  sets $G\in\cs$ with $F\prec G$ and $|G|=|F|+1$.  Thus, if $o(T_F)=\b+1$,
  then $o(T_G)\leq\b$ for every such set $G$.
  
  To simplify the notation for the induction we shall use an enumeration
  $(G_j)_0^{\infty}$ of \cs, the domain of $\psi$, which satisfies:
  $G_0=\emptyset$, $G_1=\{1\}$, and if $\max G_j<\max G_k$, then $j<k$.
  
  Now let us inductively define $\psi$.  Let
  $\psi\emptyset=\emptyset=F_0$, $\psi G_1=\psi \{1\}=F_1$, and set
  $k_1=1$.  Suppose that $\psi G_i$ has been defined for $i\leq n$ such
  that if $i<j$, $\psi G_i=F_{k_i}$ and $\psi G_j=F_{k_j}$, then
  $k_i<k_j$, $G_i\prec G_j$ if and only if $F_{k_i}\prec F_{k_j}$ and for
  $i=1,\dots,n$, $o(T_{G_i})\leq o(T_{F_{k_i}})$.  We next define $\psi
  G_{n+1}$.  From our enumeration $(G_j)_0^{\infty}$ of \cs\ there exists
  $m\geq1$ such that $G_{n+1}=\{m\}$, or $G_{n+1}=G_i\cup\{m\}$ for some
  $i\leq n$.  In the first case let $k_{n+1}>k_n$ be the least integer
  such that $|F_{k_{n+1}}|=1$ and $o(T_{F_{k_{n+1}}})\geq o(T_{G_{n+1}})$.
  We can achieve this last condition because $\sup_{r\geq1} o(T_{\{r\}}) =
  \w^{\a}$.  In the second case let $k_{n+1}>k_n$ be least with
  $F_{k_i}\prec F_{k_{n+1}}$, $|F_{k_{n+1}}|=|F_{k_i}|+1$, and
  $o(T_{F_{k_{n+1}}})\geq o(T_{G_{n+1}})$.  The existence of $k_{n+1}$ is
  guaranteed by the conditions on $\psi$.  Indeed, since
  $G_{n+1}=G_i\cup\{m\}$, it follows that $G_i$ is not a terminal node of
  $T_{\emptyset}$ so that $o(T_{G_i})>1$.  We also have
  $o(T_{F_{k_i}})\geq o(T_{G_i})$, so that neither is $F_{k_i}$ a terminal
  node of $T_{\emptyset}$ and hence $F_{k_i}$ has an infinite sequence of
  successor nodes $(E_j)_1^{\infty}$ such that $E_j\in\cs$, $F_{k_i}\prec
  E_j$, $|E_j|=|F_{k_i}|+1$ for each $j\geq1$ and $\sup_j o(T_{E_j}) =
  o(T_{F_{k_i}}) - 1 = \beta$, where $o(T_{F_{k_i}}) = \beta + 1$.  
  If \b\  is a successor ordinal, then choose $E_j$ so that
  \[ o(T_{E_j}) = \beta = o(T_{F_{k_i}}) - 1 \geq o(T_{G_i}) - 1 \geq
  o(T_{G_{n+1}}) \ . \] 
  Otherwise $o(T_{G_{n+1}}) < \beta = o(T_{G_i}) - 1$, and we may choose
  $E_j$ so that $o(T_{E_j}) \geq o(T_{G_{n+1}})$.  We then set
  $F_{k_{n+1}}=E_j$.
  
  Clearly $\psi G_n=F_{k_n}\ (n\geq1)$ satisfies all the requirements of
  the induction.  It remains to show that this is sufficient to ensure
  that the map $U$ is an isometry.  We must show that $\|\sum
  a_iU(f_i)\|_{\C} = \sup_{F\in\cs} |\sum a_i|$.  Now,
  \begin{align*}
    \Bigl\| \sum a_i U(f_i) \Bigr\|_{\C} 
      &= \Bigl\| \sum a_i \sum_{\substack{G_j\in\cs \\ \max G_j=i}} 
         \chi_{\psif{G_j}} \Bigr\|_{\C} \\  
      &= \sup_{F\in\cs} \Bigl| 
         \sum a_i\sum_{\substack{G_j\in\cs \\ \max G_j=i}}
                               \chi_{\psif{G_j}}(F) \Bigr| \ . 
   \end{align*}
   First note that for $i$ fixed, if there exist $j,j'$ such that
   $\chi_{\psif{G_j}}(F)=\chi_{\psif{G_{j'}}}(F)=1$, and $\max G_j = \max
   G_{j'} = i$, then both $\psi G_j\prq F$ and $\psi G_{j'}\prq F$ so that
   we may assume $\psi G_j\prq \psi G_{j'}$.  By the conditions on $\psi$
   this forces $G_j\prq G_{j'}$, but $\max G_j = \max G_{j'}$ so that
   $G_j=G_{j'}$, and hence $j=j'$.  Thus, for each $i\geq1$, $F\in\cs$ we
   have $\sum_{G_j\in\cs, \ \max G_j=i} \chi_{\psif{G_j}}(F)=0$ or $1$.

  To complete the proof we define a map $\varphi$ from \cs\ into the
  collection of finite subsets on \N\ by 
  \[ \varphi(F) = \Bigl\{ i\geq 1 : 
                     \sum_{\substack{G_j\in\cs \\ \max G_j=i}}
                     \chi_{\psif{G_j}}(F) = 1 \Bigr\}
  \]
  and show that the range of $\varphi{}$ is \cs, for then 
  \[ \sup_{F\in\cs} \Bigl| \sum a_i 
     \sum_{\substack{G_j\in\cs \\ \max G_j=i}}
     \chi_{\psif{G_j}}(F) \Bigr| = \sup_{F\in\cs} \Bigl|\sum a_i\Bigr| 
  \]
  as required.  

  First note that if $E\prq F$, then $\varphi(E)\prq\varphi(F)$.  Indeed, 
  \begin{align*}
    \varphi(F) &= \Bigl\{i\geq 1 :\sum_{\substack{G_j\in\cs\\\max G_j=i}} 
                     \chi_{\psif{G_j}}(F) = 1 \Bigr\}                  \\ 
               &= \{ i\geq  : \text{ there exists } j\geq1 \text{ with } i
                     = \max G_j 
                     \text{ and } \psi G_j\prq F \} \\
               &= \{ \max G_j :\psi G_j\prq F\}\ . \tag{$*$}\label{eq:phi}
  \end{align*}
  Now fix $j\geq1$ and let $E_1\prec\dots\prec E_k = G_j$ satisfy
  $|E_1|=1$ and $|E_{i+1}|=|E_i|+1\ (i<k)$.  Then $\psi E_i\prq\psi G_j$,
  so that $\max E_i\in\varphi(\psi G_j)\ (i=1,\dots,k)$, i.e.\ 
  $G_j\subseteq\varphi(\psi G_j)$.  On the other hand, if $\psi
  G_{j'}\prq\psi G_j$, then $G_{j'}\prq G_j$ which gives $G_{j'}=E_i$ for
  some $i$, so that $\varphi(\psi G_j)\subseteq G_j$ by~(\ref{eq:phi}),
  and hence the two sets are equal.  Thus \cs\ is contained in the range
  of $\varphi$.  Finally let $F\in\cs$, set $i = \max \varphi(F)$ and find
  $j_0$ such that $i = \max G_{j_0}$ and $\psi G_{j_0} \prq F$.  But then
  $G_{j_0}=\varphi(\psi G_{j_0})\prq \varphi(F)$, while $\max
  G_{j_0}=i=\max \varphi(F)$.  Hence $G_{j_0}=\varphi(F)$ and so the range
  of $\varphi$ is exactly \cs\ as required.  This completes the proof.
\end{proof}

\begin{lem}\label{lem:Xa-Ca-llbbis}
  For $\a<\w_1$, $\omega^{\alpha+1} \leq I_b(X_{\alpha}) \leq I_b(\C) \leq
  \omega^{\alpha+2}$.  Moreover, when $\omega\leq\alpha<\omega_1,\ 
  I_b(X_{\alpha})=I_b(\C)=\omega^{\alpha+1}$.
\end{lem}
\begin{proof}
  First let $(e_i)_1^{\infty}$ be the unit vector basis for $X_{\alpha}$
  and set $T=\{ (e_i)_{i\in F} : F\in{\mathcal S}_{\a} \}$.  The tree $T$
  is clearly an \llbbt\ on $X_{\alpha}$ isomorphic to ${\rm
    Tree}({\mathcal S}_{\a})\setminus\{\emptyset\}$, so that
  $o(T)=\omega^{\alpha}$.  By Lemma~\ref{lem:bigger-P-tree} the block
  basis index is strictly greater than the order of any block basis tree
  on the space, so that $\omega^{\alpha}<I_b(X_{\alpha})$.  But now, by
  Corollary~\ref{cor:I-P-bb-index}, the block basis index is of the form
  $\omega^{\beta}$ for some $\beta<\omega_1$ so that $\omega^{\a+1}\leq
  I_b(X_{\a})$.
  
  As we noted after Theorem~\ref{thm:bp}, $I^+_w(\C)=\eta(\C)=\w^{\a+1}$,
  and since the node basis for \C\ is shrinking, it follows from
  Theorem~\ref{thm:I+w<=I+b} that $I^+_b(\C)=\omega^{\alpha+1}$ when
  $\w\leq\a$ and $I^+_b(\C)\leq \omega^{\alpha+2}$ when $\a=n<\w$.  It is
  clear that $I_b(\C)\leq I^+_b(\C)$, and finally we showed in
  Lemma~\ref{lem:Xa-bb-Cwwa} that $X_{\alpha}$ embeds into \C\ as a block
  basis, and hence we have the inequalities
  \begin{gather*}
    \omega^{\a+1}\leq I_b(X_{\a})\leq I_b(\C)\leq \omega^{\a+1} 
       \text{ when $\w\leq\a$, and} \\
    \omega^{\a+1}\leq I_b(X_{\a})\leq I_b(\C)\leq \omega^{\a+2}
       \text{ when $\a=n<\w$,}
  \end{gather*}
which completes the proof.
\end{proof}
\begin{rem}\label{rem:Xa-Ca-llis}
  Since $\omega\leq\alpha$, it follows from
  Theorem~\ref{thm:j-reflexive}~(\ref{thm:j-reflexive:item:B}) that also
  $I(X_{\alpha})=I(\C)=\omega^{\alpha+1}$.
\end{rem}

\begin{lem}\label{lem:Cwwn++Cwwnk->Cwwn}
  For each $n\geq1$, every $k\geq1$ and any admissible enumerations of the
  node bases of $(\C[n]\oplus\dots\oplus\CC[n]{k})_{\infty}$ and \C[n],
  the node basis of $(\C[n]\oplus\dots\oplus\CC[n]{k})_{\infty}$ embeds
  isomorphically into $\C[n]$ as a block basis of the node basis.
\end{lem}

Before we can prove this lemma we need to extend the definition of node
basis from \C\ to \CC{k}\ and $(\C\oplus\dots\oplus\CC{k})_{\infty}$ for
$\a<\w,\ k\geq1$.  First observe that in \cs[\a+1]\ we have a natural copy
of $\w^{\w^{\a}\cdot k}$ given by 
\[ \cs[\a,k] = \{ \{k+1\}\cup\bigcup_{i=1}^k F_i : 
   \{k+1\}<F_1<\dots<F_k,\    F_i\in\cs \} \]
so that $\{ \ch_{\{G:F\prq G\}} = \chir_F : F\in\cs[\a,k] \}$ is
the node
basis for \CC{k}.  Further, $\cs[\a,k]\cap\cs[\a,l]=\emptyset$ if $k\neq
l$, thus $\{ \chir_{\{l\}} \}_{l=2}^{k+1}$ is a sequence of disjointly
supported 
functionals and $\{ \chir_F : \{l\}\prq F, 2\leq l\leq k+1, F\in\cs[\a+1]
\}$ is a node basis for 
\[
  (\C\oplus\dots\oplus\CC{k})_{\infty} 
     = \{ f\in\C[\a+1] : f(F)=0 \text{ if there exists } j>k+1
          \text{ with } \{j\}\prq F \} \ . 
\]
The natural projection $Q_k$ of \C[\a+1]\ onto
$(\C\oplus\dots\oplus\CC{k})_{\infty}$ is given by
\[ Q_k g = \Bigl( \sum_{l=2}^{k+1} \chir_{\{l\}} \Bigr)\cdot g \ . \]
Finally we note that if $(e_i)_0^{\infty}$ is an admissible ordering of
the node basis of \C, then $e_0=\chir_{\emptyset}=\ch_{[1,\w^{\w^{\a}}]}$,
and hence $(e_i)_1^{\infty}$ is a node basis for $\Co=\{f\in\C :
f(\w^{\w^{\a}}) = 0 \}$.

\begin{proof}[Proof of Lemma~\ref{lem:Cwwn++Cwwnk->Cwwn}]
  The argument follows the same lines as the proof that \CC{k}\ is
  isomorphic to \C\ in~\cite{bp}.  Note that for $\a=0$ the node basis of
  $(C(\w)\oplus\dots\oplus C(\w^k))_{\infty}$ is a family of indicator
  functions with nested or disjoint supports, and the nested functions are
  at most $k+1$ sets deep.  The required map $T$ is found by sending the
  $i^{\text{th}}$ element of the admissible enumeration of the node basis
  of $(C(\w)\oplus\dots\oplus C(\w^k))_{\infty}$ to the
  $(i+1)^{\text{th}}$ element $\chir_{\{i\}}$ of the node basis of
  $C(\w)$, $(\chir_{\emptyset},\chir_{\{1\}},\chir_{\{2\}},\dots)$.  (Note
  that $\chir_{\emptyset}$ is not in the image.)  It is easy to see that
  $\|T\|\cdot\|T^{-1}\|\leq2(k+1)$.  The general case is similar.
  
  We view the node basis of $(\C[n]\oplus\dots\oplus\CC[n]{k})_{\infty}$
  (in the ordering $\prec$) as a disjoint union of $k$ trees
  $\{\chir_{\{m+1\}}\}\cup T(m)$, $m=1,2,\dots,k$, with $T(m)$ isomorphic
  to the replacement tree $T(m,\w^{\w^n})$ and $\chir_{\{m+1\}}$ the
  unique initial node of the tree $\{\chir_{\{m+1\}}\}\cup T(m)$.  Thus
  $z\in T(m)$ implies $z=\chir_F$ with $\{m+1\}\prec F$.
  
  For each $m=1,\dots,k$ let $F_m : T(m)\rightarrow
  T_m=\{a_1^m,\dots,a_m^m\}$ be the defining map for the replacement tree.
  Recall that $F^{-1}_m(a_i^m)$ is one or a countable union of trees, each
  isomorphic to $T_{\w^{\w^n}}$.  Let $(U_j)_{j=1}^{\infty}$ be an
  enumeration of all of these trees for $1\leq i\leq m\leq k$.  For each
  $j$ let $(\chir_{G_{jl}})_{l=1}^{\infty}$ be the sequence of initial
  nodes of $U_j$, so that $\{ \chir_F\in U_j : G_{jl}\prq F \}$ is
  equivalent to the node basis of \CC[n-1]{l}\ under the natural map.  Let
  $(y_i)_1^{\infty}$ be the given admissible enumeration of the node basis
  of $(\C[n]\oplus\dots\oplus\CC[n]{k})_{\infty}$ and let
  $(w_j)_1^{\infty}$ be an admissible enumeration of the node basis of
  \C[n].  To avoid confusion between domain and range we shall let
  $\zetar_F = \ch_{\{G\in\cs[n] : F\prq G \}}$, for $F\in\cs[n]$ denote
  the elements of the node basis of \C[n]\ in the image.  Thus $\{w_j :
  j\geq1 \} = \{ \zetar_F : F\in\cs[n] \}$.

  We define a map $\psi:(y_i)_1^{\infty}\rightarrow(w_j)_1^{\infty}$
  inductively to satisfy the following conditions:
  \renewcommand{\labelenumi}{(\roman{enumi})}
  \begin{enumerate}
  \item $w_1=\zetar_{\emptyset}$ is not in the image of $\psi$;
  \item if $y_i=\chir_E$ and $E=\{m\}$ $(m=1,\dots,k)$ or $E=G_{jl}$ for
    some $j,l\geq1$, then $\psi(y_i)=\zetar_{\{s\}}$ for some $s\geq1$;
  \item $\psi$ is increasing, i.e.\ if $\psi(y_i)=w_{l(i)}$, then
    $l(1)<l(2)<\cdots$; 
  \item if $G_{jl}\prec E_1\prec E_2$ and
    $\psi(\chir_{G_{jl}})=\zetar_{F_0},\ \psi(E_i)=\zetar_{F_i}\ (i=1,2)$,
    then $F_0\prec F_1\prec F_2$;
  \item if $\psi(E)=\zetar_F,\ G_{jl}\prec E$ and $\{s\}\prec F$, then the
    order of $\chir_E$ in $\{\chir_H : G_{jl}\prq H, H\in U_j\}$ is less
    than or equal to the order of $\zetar_F$ in $\{\zetar_H : \{s\}\prq
    H\}$, where the sets are trees in the usual order $\prec$ and the
    order of a node $z$ in a tree $T$ is simply the order of the subtree
    $\{y\in T : y\leq z\}$ of $T$.
  \end{enumerate}
  It is easy to see that the inductive definition of $\psi$ will succeed
  because if $\psi(y_1),\dots,\psi(y_i)$ have been chosen, then there are
  infinitely many candidates for $\psi(y_{i+1})$ satisfying (i)--(v).  It
  is also not difficult to see that if $S$ is the induced map from 
  $(\C[n]\oplus\dots\oplus\CC[n]{k})_{\infty}$ into \C[n], then
  $\|S\|\cdot\|S^{-1}\|\leq2(k+1)$. 
\end{proof}

\begin{rem}\label{rem:Cwwn++Cwwnk->Cwwn}
  It is clear from the proof that the blocking of the basis of \C[n]\ is
  actually just a subsequence.  The same argument works for
  $(\C\oplus\dots\oplus\CC{k})_{\infty}$ into \C, and the argument also
  shows that the node basis of \C[n]\ is equivalent to a subsequence of
  the node basis of \Co[n].
\end{rem}

\begin{lem}\label{lem:Xn-Cn-llbbi}
  For $n\geq1,\ I_b(X_n)=I_b(\C[n])=\omega^{n+1}$.
\end{lem}
\begin{proof}
  By Lemma~\ref{lem:Xa-bb-Cwwa} and the proof of
  Lemma~\ref{lem:Xa-Ca-llbbis} we have
  \[ \omega^{n+1} \leq I_b(X_n) \leq I_b(\C[n]) \ . \]
  To complete the proof we show that for each $n\geq0$ there does not
  exist an \llbbt\ on $\Co[n]$ of order $\omega^{n+1}$, and hence
  $I_b(\Co[n]) \leq \omega^{n+1}$.  Then, since
  $I_b(X,(e_i)_0^{\infty})=I_b(X,\seq{e_i})$ for any space $X$ with basis
  $(e_i)_0^{\infty}$, and $\Co[n]=[e_i]_1^{\infty}$ where
  $(e_i)_0^{\infty}$ is any admissible enumeration of the node basis for
  \C[n], it follows that $I_b(\C[n])=I_b(\Co[n])\leq\w^{n+1}$.
  
  We prove this result by induction on $n$.  For $n=0$ we first note that
  $C_0(\omega)=c_0$.  Since the unit vector basis of $c_0$ does not
  contain $\ell_1^n$'s uniformly as block bases, it follows that $c_0$
  contains no \llbbt\ of order $\omega$.

  We assume that the result is
  true for $n$, and let $\{e_i : i\geq1\}$ be an admissible enumeration of
  the node basis of \Co[n+1].  Suppose that $T$ is an \llkbbt\ of order
  $\w^{n+2}$ on \Co[n+1]\ which, without loss of generality, we assume
  consists of finitely supported vectors with respect to \seq{e_i}, and is
  isomorphic to the minimal replacement tree $T(\w,\w^{n+1})$.

  We write $T=\cup_{m=1}^{\infty}T(m)$, where $T(m)$ is a tree isomorphic
  to $T(m,\w^{n+1})$ and the elements from different trees $T(m)$ are
  unrelated.  Choose $m>2K$ and let $F:T(m)\rightarrow
  T_m=\{a_1,\dots,a_m\}$, where $a_1<a_2<\dots<a_m$, be the defining map
  for the replacement tree $T(m,\w^{n+1})$.  Let $S_1=F^{-1}(a_1)$ so that
  $S_1\simeq T_{\w^{n+1}}$.  Let $(x_i^1)_{i=1}^{p_1}$ be a terminal node
  in $S_1$.  Define $x_1=x_1^1$ and let 
  \[ k_1 = \max \{ k\geq1 : Q_k x_1 \neq 0 \} \ . \]
  
  Let $S'_1 = \{ z\in T(m) : (x_i^1)_{i=1}^{p_1} < z \}$, so that
  $S'_1\cap F^{-1}(a_2)$ is isomorphic to $T_{\w^{n+1}}$, and let $S_2$ be
  the restricted tree $R(S'_1\cap F^{-1}(a_2))$.  The tree $S_2$ is an
  \llbbt\ of order $\w^{n+1}$.  By Lemma~\ref{lem:Cwwn++Cwwnk->Cwwn} and
  the induction hypothesis there is no \llbbt\ on $Q_{k_1}(\Co[n+1])$ of
  order $\w^{n+1}$.  Consider the tree
  \[ Q_{k_1} S_2 = \{ (Q_{k_1} z_1,\dots,Q_{k_1} z_l) : (z_1,\dots,z_l)\in
  S_2 \} \ ; \] note that $\|Q_{k_1} z_i\|\leq\|Q_{k_1}\| \|z_i\|\leq1$
  and $\spp Q_{k_1} z_i \cap Q_{k_1} z_j = \emptyset$ when $i\neq j$ and
  $(z_i)_1^l\in S_2$ (since $S_2$ is a block basis tree).  If there exists
  $\delta>0$ such that for every $(z_i)_1^l\in S_2$ and
  $(b_i)_1^l\subset\R$,
  \[ \Bigl\| \sum_{i=1}^l a_i Q_{k_1} z_i \Bigr\| 
             \geq \delta \sum_1^l |b_i| \ , \]
  then the tree
  \[ \{ (Q_{k_1} z_1/\| Q_{k_1} z_1 
  \|,\dots,Q_{k_1} z_l/\| Q_{k_1} z_l \|) : (z_i)_1^l\in S_2 \} \] 
  would be an $\ell_1$-$\delta^{-1}$-block basis tree on $Q_{k_1}
  \Co[n+1]$ of order $\w^{n+1}$, contradicting the induction hypothesis.
  Therefore there is a terminal node $(x_i^2)_{i=1}^{p_2}\in S_2$ and
  $(b_i^2)_{i=1}^{p_2}\subset\R$ such that $\sum_1^{p_2}|b_i^2|=1$ and
  $\|\sum_1^{p_2} b_i^2 Q_{k_1} x_i^2 \|< 1/m$.  Define $x_2 =
  \sum_1^{p_2} b_i^2 x_i^2$ and let $k_2 = \max \{ k\geq1 : Q_k x_2 \neq 0
  \} \vee (k_1+1)$.
  
  As before we consider the tree $S'_2 = \{ z\in T(m) :
  (x_1^1,\dots,x_{p_1}^1,x_1^2,\dots,x_{p_2}^2) < z \}$, so that $S'_2\cap
  F^{-1}(a_3)$ is isomorphic to $T_{\w^{n+1}}$, and we let $S_3$ be the
  restricted tree $R(S'_2\cap F^{-1}(a_3))$.  Arguing as above there is a
  terminal node $(x_i^3)_{i=1}^{p_3}\in S_3$ and
  $(b_i^3)_{i=1}^{p_3}\subset\R$ such that $\sum_1^{p_3}|b_i^3|=1$ and
  setting $x_3 = \sum_1^{p_3} b_i^3 x_i^3$ gives $\| Q_{k_2} x_3 \|<1/m$.
  Continuing in this way we get $(x_i)_1^m$ a block basis of some node
  $y=(x_1^1,\dots,x_{p_1}^1, x_1^2,\dots,x_{p_2}^2, \dots,
  x_1^m,\dots,x_{p_m}^m)$ of $T(m)$ such that $x_i = \sum_1^{p_i} b_j^i
  x_j^i$ for some sequence $(b_j^i)_{j=1}^{p_i}\subset\R$ with
  $\sum_1^{p_i}|b_j^i|=1$, together with a sequence $(k_i)_1^m$ such that
  $k_i = \max \{ k\geq1 : Q_k x_i \neq 0 \} \vee (k_{i-1}+1)$ (where
  $k_0=0$).

  Let $x=\frc{m}\sum_{i=1}^m x_i$, so that
  \[
    \| x \| = \Bigl\| \frc{m}\sum_{i=1}^m x_i \Bigr\| 
            = \Bigl\| \frc{m}\sum_{i=1}^m 
                             \sum_{j=1}^{p_i} b_j^i x_j^i \Bigr\| 
            \geq \frc{K}\frc{m}\sum_{i=1}^m \sum_{j=1}^{p_i} |b_j^i|
            \geq  \frc{K}\frc{m}\sum_1^m 1 = \frc{K} \ ,
  \]
  since $T$ was an \llkt.  On the other hand (with $Q_{k_0}=0$)
  \begin{align*}
   \| x \| &= \max_{ 1 \leq j \leq m } \| (Q_{k_j} - Q_{k_{j-1}}) x \|
            = \max_{ 1 \leq j \leq m } \Bigl\| (Q_{k_j} - Q_{k_{j-1}}) 
                                     \frc{m} \sum_{i=1}^m x_i \Bigr\| \\
           &\leq \max_{ 1 \leq j \leq m } \frc{m} \Bigl( \|x_j\| + 
                          \sum_{i=j+1}^m \| Q_{k_{i-1}} x_i \| \Bigr) 
            \leq \max_{ 1 \leq j \leq m } \frc{m} 
                                      \Bigl( 1 + \frac{m-j}{m} \Bigr) 
            < \frac{2}{m} 
            < \frc{K} \ .
  \end{align*}
  Thus there exists no such tree $T$ of order $\w^{n+2}$ on \Co[n+1]\ 
  which completes the proof.
\end{proof}

The goal of the next few results is to show that the \lli\ of \C[n]\ is
$\omega^{n+2}$.  First we need some preliminary results.

\begin{lem}\label{lem:CK->CKx1..n}
  Let $K$ be a compact Hausdorff space, let $\kappa:C(K)\hookrightarrow
  C(K\times\{1,\dots,2^n\})$ be the map $(\kappa f)(k,j) = f(k)$, let
  $\io: C(\{1,\dots,2^n\})\hookrightarrow C(K\times\{1,\dots,2^n\})$ be
  the map $(\io f)(k,j) = f(j)$, and let $(r_m)_1^n$ be the standard
  Rademacher functions on $\{1,\dots,2^n\}$ so that $(r_m)_1^n$ is
  1-equivalent to the unit vector basis of $\llo^n$.  Then for every $f\in
  C(K)$ and any sequence $(a_m)_1^n\subset\R$,
  \[ \Bigl\| \kappa f + \io \sum_1^n a_m r_m \Bigr\| 
             = \| f \| + \sum_1^n | a_m | \ . \]
\end{lem}

\begin{proof}
  Find $k_0$ so that $|f(k_0)| = \|f\|$,  let
  $\e=\sgn(f(k_0))$, and find $l_0\in\{1,\dots,2^n\}$ such that $\sum_1^n
  a_m r_m(l_0) = \e \sum_1^n | a_m |$.  Now,
  \begin{align*}
    \Bigl| (\kappa f)(k_0,l_0) + 
           \Bigl( \io \sum_1^n a_m r_m \Bigr)(k_0,l_0) \Bigr|
           &= \Bigl| f(k_0) + \sum_1^n a_m r_m(l_0) \Bigr| \\
           &= \Bigl| \e \| f \| + \e \sum_1^n | a_m | \Bigr| \\
           &= \| f \| + \sum_1^n | a_m |,
  \end{align*}
  as required.
\end{proof}

Note that this result will also apply if we replace $C(K)$ with
$C_0(\wa)$, since functions on $C_0(\wa)$ attain their norm.

\begin{lem}\label{lem:incr-l1-ind}
  Let $1\leq\gamma<\alpha<\omega_1$; if there exists an \llt\ with
  constant 1 and order \b\ on $C_0(\w^{\a})$, then there exists an \llt\ 
  on $C_0(\w^{\a+\gamma})$ with constant 1 and order $\b+\ab{\w}{\gamma}$.
\end{lem}
\begin{proof}
  We may write $C_0(\w^{\a+1})$ as 
  \[ C_0(\w^{\a+1})
       = \Bigl( \sum_{j=1}^{\infty}\oplus
                C_0(\w^{\a} \times (2^j,2^{j+1}]) \Bigr)_{c_0}
       = \Bigl( C_0(\w^{\a} \times (2,4] ) \oplus 
                C_0(\w^{\a} \times (4,8] ) \oplus \dots \Bigr)_{c_0} \ , 
  \] 
  where $(2^j,2^{j+1}] = \{ 2^j+1,2^j+2,\dots,2^{j+1}\}$.  We shall prove
  the result using induction on $\gamma$.
  
  Let $\gamma=1$ and let $\kappa_j:C_0(\w^{\a})\hookrightarrow
  C_0(\w^{\a}\times(2^j,2^{j+1}])$ be the map $(\kappa_j f)(\b,l)=f(\b)$,
  restricted from $C(\wa)$ to $C_0(\wa)$ $(j=1,2,\dots)$, let $(\bar
  r^j_i)_{i=1}^j$ be the Rademacher functions on $(2^j,2^{j+1}]$ and let
  $(r^j_i)_{i=1}^j$ be the extension of these to
  $C_0(\w^{\a}\times(2^j,2^{j+1}])$ with $r^j_i(\b,l)=\bar r^j_i(l)$.  Let
  $T$ be a tree with constant 1 and order \b\ on $C_0(\w^{\a})$; we
  construct a tree of order $\b+\w$ on $C_0(\w^{\a+1})$.  Let
  \[ S_j = \{ (r_1^j),(r_1^j,r_2^j),\dots,(r_1^j,\dots,r_j^j) \} \cup
           \{ (r_1^j,\dots,r_j^j,\kappa_j x_1,\dots,\kappa_j x_m) :
              (x_i)_1^m\in T \} \ , 
  \]
  and let $S=\cup_{j=1}^{\infty} S_j$ with the usual ordering by
  extension.  The subtree of $S$ given by $\cup_{j=1}^{\infty}
  \{ (r_1^j),(r_1^j,r_2^j),\dots,(r_i^j)_{i=1}^j \}$ has order
  \w, and after every terminal node is a tree of order \b\ so that
  $o(S)=\b+\w$.  It is clear from the previous lemma that $S$ is an \llt\
  with constant 1.  

  If the result is true for $\gamma$, then given an \llo-1-tree on
  $C_0(\w^{\a})$ of order \b, there exists an \llo-1-tree on
  $C_0(\w^{\a+\gamma})$ of order $\b+\ab{\w}{\gamma}$, but now by the case
  $\gamma=1$ there exists an \llo-1-tree on
  $C_0(\w^{\a+(\gamma+1)})=C_0(\w^{(\a+\gamma)+1})$ of order
  $(\b+\ab{\w}{\gamma})+\w = \b + \ab{\w}{(\gamma+1)}$.  

  Finally, if $\gamma$ is a limit ordinal and the result has been proven
  for every $\gamma'<\gamma$, then let $\gamma_n\nearrow\gamma$ and 
  \[ C_0(\w^{\a+\gamma}) \simeq 
       \Bigl( C_0(\w^{\a+\gamma_1}) \oplus 
              C_0(\w^{\a+\gamma_2}) \oplus \dots \Bigr)_{c_0}\ , 
  \] 
  hence we may take the union of \llo-1-trees $S_n$ on
  $C_0(\w^{\a+\gamma_n})$ of order $\b+\ab{\w}{\gamma_n}$ to obtain a tree
  on $C_0(\w^{\a+\gamma})$ of order
  $\sup_n(\b+\ab{\w}{\gamma_n})=\b+\ab{\w}{\gamma}$ as required.
\end{proof}

\begin{lem}\label{lem:Cn-lli}
  $I(\C[n])=\omega^{n+2}$.
\end{lem}
\begin{proof}
  {} From Theorem~\ref{thm:j-reflexive}~(\ref{thm:j-reflexive:item:C}) and
  Lemma~\ref{lem:Xn-Cn-llbbi} we know that the \lli\ of \C[n]\ is either
  $\omega^{n+1}$ or $\omega^{n+2}$ and hence $I(\Co[n])=\omega^{n+1}$ or
  $\omega^{n+2}$.  For each $n\geq0$ we shall construct an \llt\ on
  \Co[n]\ of order $\omega^{n+1}$ so that $I(\Co[n])=\omega^{n+2}$ by
  Lemma~\ref{lem:bigger-P-tree} and the result follows.  This is clear for
  $n=0$ since $\llo^n$ embeds isometrically into $\ell^{2^n}_{\infty}$ for
  each $n\geq1$, which immediately yields an $\ell_1$-1-tree of order
  $\omega$.

  We may now complete the proof by induction on $n$.  If there is an
  \llo-1-tree on \Co[n]\ of order $\w^{n+1}$, then by the previous lemma
  there exists a tree of order $\ab{\w^{n+1}}{k} = \w^{n+1} +
  \ab{\w}{(\ab{\w^n}{(k-1)})}$ on $\CCo[n]{k} =
  C_0(\w^{\w^n+\w^n\cdot(k-1)})$ for every $k\geq1$.  Taking the union
  over $k$ of these we obtain an \llo-1-tree on \Co[n+1]\ of order
  $\w^{n+2}$ as required.  This completes the inductive step and hence the
  proof.
  \end{proof}

\begin{lem}\label{lem:Xn-lli}
  $I(X_n)=\omega^{n+1}$.
\end{lem}
\begin{proof}
  Again, from Theorem~\ref{thm:j-reflexive}~(\ref{thm:j-reflexive:item:I})
  and Lemma~\ref{lem:Xn-Cn-llbbi} we know that $I(X_n)$ is either
  $\omega^{n+1}$ or $\omega^{n+2}$.  To demonstrate that it is the former
  we show that for each $n\geq1$ there does not exist an \llt\ on $X_n$ of
  order $\omega^{n+1}$.
 
  We prove this by induction on $n$ based on the following lemmas.  The
  idea of the proof is that if we do have an \llt\ of order $\omega^{n+1}$
  on $X_n$, then we can find a node in that tree which admits an absolute
  convex combination with arbitrarily small norm.  This contradicts the
  hypothesis that it was an \llt.
\end{proof}

Below, if $x\in X_1=[e_i]$, with $x=\sum a_ie_i$, then we define the
supremum norm of $x$ to be $\|x\|_{\infty} = \sup|a_i|$.

\begin{lem}\label{lem:inf-norm<e}
  For each $\e>0$ and each $K\geq1$ there exists $n\geq1$ such that if
  $(x_i)_1^n$ is a sequence of norm one vectors in $X_1$ which is
  $K$-equivalent to the unit vector basis of $\ell_1^n$, then there exists
  a norm one vector $x\in S([x_i]_1^n)$ with $\|x\|_{\infty}<\e$.
\end{lem}
\begin{proof}
  Fix $n$ and let $(x_i)_1^n$ be as in the statement of the lemma.
  Suppose that $\|x\|_{\infty}\geq\e$ for each $x\in S([x_i]_1^n)$.  Then
  $\|x\|_{X_1}\geq\|x\|_{\infty}\geq\e\|x\|_{X_1}$ for each $x\in [x_i]$.
  We may assume that each $x_i$ has finite support with respect to the
  unit vector basis of $X_1$, and let $N=\max\{ \spp(x_i) : i\leq n \}$.
  Thus $([x_i]_1^n,\|\cdot\|_{X_1})$ embeds into $\ell_{\infty}^N$ with
  constant $1/\e$ via the map $\ \hat{ } : x\mapsto(e^*_j(x))_{j=1}^N$,
  and hence $(\hat{x}_i)_1^n$ has a lower $\ell_1$ estimate with constant
  $\e/K$.
  
  By James \cite{j}, for fixed $k$ and $\delta>0$, if $n$ is sufficiently
  large, then there exists a normalized block basis $(\hat{y}_i)_1^k$ of
  $(\hat{x}_i)_1^n$ such that
  $(\hat{y}_i)_1^k\stackrel{1+\delta}{\sim}\mbox{uvb }\ell_1^k$.  Now if
  we take $\delta$ to be very small, depending on $k$, then we see that
  for each $i=1,\dots,k$ the size of one of the sets $E_i=\{ j\leq N :
  \hat{y_i}(j)>1/2 \}$ and $F_i=\{ j\leq N : \hat{y_i}(j)<-1/2 \}$ must be
  at least $2^{k-2}$.  We calculate the norm of $y_1$ in $X_1$ supposing
  that $|E_1|\geq2^{k-2}$.  Let $E$ be the second half of $E_1$, so that
  if $E_1=\{e_1,\dots,e_r\}$, then $E=\{e_s,e_{s+1},\dots,e_r\}$, where
  $s=[(r+1)/2]$.  Clearly $E\in\cs[1]$,
  $|E|\geq\frac{1}{2}2^{k-2}=2^{k-3}$ and \[\|y_1\|_{X_1}\geq
  \Bigl|\sum_{j\in E}y_1(j)\Bigr|=\Bigl|\sum_{j\in E}\hat{y}_1(j)\Bigr| 
  \geq|E|\!\cdot\!\frac{1}{2}  \geq\frac{1}{2}2^{k-3}=2^{k-4}\ .\]
  On the other hand
  $\|y_1\|_{X_1}\leq\frac{1}{\e}\|\hat{y}_1\|_{\infty}=\frac{1}{\e}$, so
  this is impossible for large $k$, and hence for $n$ large enough.

  This contradicts our initial assumption that $\|x\|_{\infty}\geq\e$ for
  each $x\in S([x_i]_1^n)$ and hence there exists $x\in S([x_i]_1^n)$ with
  $\|x\|_{\infty}<\varepsilon$.  
\end{proof}

\begin{lem}\label{lem:x1-inf-norm<e}
  If $T$ is a tree on $B_{X_1}\setminus\{0\}$ of order $\omega$, then for
  any $\e>0$ there exist $(x_i)_1^n\in T$ and $(a_i)_1^n\subset\R$ with
  $\sum|a_i|=1$ and $\|\sum_1^na_ix_i\|_{\infty}<\e$.
\end{lem}
\begin{proof}
  Choose $\e>0$, set $K=1/\e$, and let $T$ be a tree on $X_1$ as above.
  If there exist $(x_i)_1^n\in T$ and $(a_i)_1^n\subset\R$ such that
  \[ \Bigl\|\sum_1^na_ix_i\Bigr\|<
  \frac{1}{K}\sum_1^n|a_i|\!\cdot\!\|x_i\|\ , \] 
  then set $a=\sum_1^n|a_i|$ and $\bar a_i=a_i/a$, and we have 
  \[
    \Bigl\|\sum_1^n\bar a_ix_i\Bigr\|_{\infty}  \leq \Bigl\|\sum_1^n\bar
                      a_ix_i\Bigr\| 
              < \frac{1}{K}\sum_1^n\frac{|a_i|}{a}\!\cdot\!\|x_i\|
              \leq \frac{1}{K}\sum_1^n\frac{|a_i|}{a} = \frac{1}{K}=\e\ ,
  \]
  while $\sum_1^n|\bar a_i|=1$ as required.
  
  Otherwise set $\bar T=\{ (\bar x_i)_1^n : (x_i)_1^n\in T \}$, where
  $\bar{x}=x/\|x\|$.  Then $\bar T$ is an \llkt\ on $X$ of order $\w$, and
  {} from the previous lemma can find $(\bar x_i)_1^n\in T$ and $
  x=\sum_1^n a_i\bar x_i\in S([\bar x_i])$ such that $\|
  x\|_{\infty}<\e$.  Now, $\| x\|=1$ and $\| x\|\leq\sum|a_i|$ so
  that $\sum_1^n|a_i|\geq1$; also $\|x_i\|\leq1$ for each $i$, thus if we
  set $a=\sum_1^n(|a_i|/\|x_i\|)$, then $a\geq1$.  Clearly
  $\frac{1}{a}\sum_1^na_i\bar x_i$ has supremum norm less than that of $x$
  so that $\|\frac{1}{a}\sum_1^na_i\bar x_i\|_{\infty}<\e$.  Finally
  \[ \frac{1}{a}\sum_1^na_i\bar x_i=\sum_1^n\frac{a_i}{a\|x_i\|}x_i\ .\] 
  so that setting $b_i=a_i/(a\|x_i\|)$, we obtain
  $\frac{1}{a}\sum_1^na_i\bar x_i=\sum b_ix_i$ with $\|\sum b_ix_i\|<\e$
  and $\sum|b_i|=1$ as required.
\end{proof}

\begin{lem}\label{lem:xn-norm<e}
  If $T$ is a tree on $B_{X_n}\setminus\{0\}$ of order $\omega^{n+1}$, and
  $\e>0$, then there exist $(x_i)_1^n\in T$ and $(a_i)_1^n\subset\R$ with
  $\sum_1^n|a_i|=1$ and $\|\sum_1^na_ix_i\|_{X_n}<\e$.
\end{lem}
\begin{proof}

  We shall prove the result by induction on $n$.  Let $n=1$ and let $T$ be
  a tree on $X_1$ of order $\omega^2$ satisfying the hypotheses of the
  lemma.  We may assume that $T$ consists of finitely supported vectors
  with respect to the basis of $X_1$, and that $T$ is isomorphic to the
  minimal tree $T(\omega,\omega)$.  We write $T=\cup_mT(m)$ where $T(m)$
  is a tree isomorphic to $T(m,\omega)$ and the elements from different
  trees $T(m)$ are unrelated.
  
  Choose $m>2/\e$, and let $F:T(m)\rightarrow T_m=\{ a_1,\dots,a_m \}$,
  where $a_1<\dots<a_m$, be the defining map for the replacement tree
  $T(m,\omega)$.  Let $S_1=F^{-1}(a_1)$, so that $S_1\simeq T_{\omega}$.
  Let $(x^1_i)_1^{p_1}$ be any terminal node in $S_1$.  Define $x_1=x^1_1$
  and let $k_1=\max(\spp x_1)$.  
  
  Let $S'_1 = \{ z\in T(m) : (x_i^1)_{i=1}^{p_1} < z \}$, so that
  $S'_1\cap F^{-1}(a_2)$ is isomorphic to $T_{\w}$, and let $S_2$ be the
  restricted tree $R(S'_1\cap F^{-1}(a_2))$.  The tree $S_2$ on $X_1$ has
  order \w, and satisfies the conditions of the previous lemma, thus we
  may find a terminal node $(x_i^2)_{i=1}^{p_2}\in S_2$ and
  $(a_i^2)_{i=1}^{p_2}\subset\R$ such that $\sum_1^{p_2}|a_i^2|=1$ and
  $\|\sum_1^{p_2} a_i^2 x_i^2 \|_{\infty}< 1/(2k_1)$.  Define $x_2 =
  \sum_1^{p_2} a_i^2 x_i^2$ and let $k_2 = \max (\spp x_2) \vee (k_1+1)$.
  
  As before we consider the tree $S'_2 = \{ z\in T(m) :
  (x_1^1,\dots,x_{p_1}^1,x_1^2,\dots,x_{p_2}^2) < z \}$, so that $S'_2\cap
  F^{-1}(a_3)$ is isomorphic to $T_{\w}$, and we let $S_3$ be the
  restricted tree $R(S'_2\cap F^{-1}(a_3))$.  Arguing as above there is a
  terminal node $(x_i^3)_{i=1}^{p_3}\in S_3$ and
  $(a_i^3)_{i=1}^{p_3}\subset\R$ such that $\sum_1^{p_3}|a_i^3|=1$ and
  setting $x_3 = \sum_1^{p_3} a_i^3 x_i^3$ gives $\| x_3 \|<1/(4k_2)$.

  Continuing in this way we get $(x_i)_1^m$ a block basis of some node
  $z=(x_1^1,\dots,x_{p_1}^1, \dots,
  x_1^m,\dots,x_{p_m}^m)$ of $T(m)$ such that $x_i = \sum_1^{p_i} a_j^i
  x_j^i$ for some sequence $(a_j^i)_{j=1}^{p_i}\subset\R$ with
  $\sum_1^{p_i}|a_j^i|=1$, together with a sequence $(k_i)_1^m$ such that
  $k_i = \max (\spp x_i) \vee (k_{i-1}+1)$, where $k_0=0$, and
  $\|x_{i+1}\|_{\infty}<1/(2^ik_i)\ (i>1)$.  Let $x=\frc{m}\sum_{i=1}^m
  x_i$, so that for $E\in\cs[1]$, if $k=\min E\geq|E|$, and $i\leq m$ is
  chosen so that $k_{i-1}<k\leq k_i$, then
  \begin{align*}
    \Bigl| \sum_{r\in E} e^*_r(x) \Bigr| 
           &\leq \frc{m} \Bigl( \|x_i\| 
                 + |E| \sum_{j=i+1}^m \|x_j\|_{\infty} \Bigr) \\
           &\leq \frc{m} \Bigl( 1 
                 + k_i \sum_{j=i+1}^m \frc{2^{j-1}k_{j-1}} \Bigr) \\
           &\leq \frc{m} \Bigl( 1 
                 + k_i \sum_{j=i+1}^m \frc{2^{j-1}k_i} \Bigr) \\
           &\leq \frac{2}{m} \ .
  \end{align*}
  Thus $\|x\|_{X_1}\leq2/m<\e$ so that $x$ is the vector we seek.  This
  completes the proof in the case $n=1$.
  
  We next suppose the result has been proven for $n-1$ and let $T$ be a
  tree on $X_n$ of order $\w^{n+1}$.  We may assume consists of
  finitely supported vectors in $X_n$, is isomorphic to
  $T(\omega,\omega^n)$, and may be written as $\cup_mT(m)$ with
  $T(m)\simeq T(m,\omega^n)$.  As before, choose $m>2/\e$, and let
  $F:T(m)\rightarrow T_m=\{ a_1,\dots,a_m \}$ be the defining map for the
  replacement tree $T(m,\omega^n)$.  Define $S_1, (x_i^1)_{i=1}^{p_1},
  x_1, k_1, S'_1$ as for the case $n=1$.  This time the tree
  $S_2=R(S'_1\cap F^{-1}(a_2))$ has order $\w^n$ on $X_n$, but we may also
  consider it as a tree of order $\w^n$ on $B_{X_{n-1}}\setminus\{0\}$, and
  hence it satisfies the conditions of the lemma for $n-1$.  Thus, by the
  induction hypothesis, there exists a terminal node
  $(x_i^2)_{i=1}^{p_2}\in S_2$ and $(a_i^2)_{i=1}^{p_2}\subset\R$ such
  that $\sum_1^{p_2}|a_i^2|=1$ and $\|\sum_1^{p_2} a_i^2 x_i^2
  \|_{X_{n-1}}< 1/(2k_1)$.

  Continuing in this way we obtain $(x_i)_1^m$ a block basis of a node
  $z=(x_1^1,\dots,x_{p_1}^1, \dots,
  x_1^m,\dots,x_{p_m}^m)$ of $T(m)$ such that $x_i = \sum_1^{p_i} a_j^i
  x_j^i$ for some sequence $(a_j^i)_{j=1}^{p_i}\subset\R$ with
  $\sum_1^{p_i}|a_j^i|=1$, together with a sequence $(k_i)_0^m$ such that
  $k_0=0$, $k_i = \max (\spp x_i) \vee (k_{i-1}+1)$ and
  $\|x_{i+1}\|_{X_{n-1}}<1/(2^ik_i)$.  

  Let $x=\frc{m}\sum_{i=1}^m
  x_i$, let $E\in\cs[n]$, $k=\min E$ and choose $i\leq m$ so that
  $k_{i-1}<k\leq k_i$.  We may write $E=\cup_{l=1}^k E_l$ where
  $E_l\in\cs[n-1]\ (l=1,\dots,k)$ and $k\leq E_1<\dots<E_k$.  Now,
  \begin{align*}
    \Bigl| \sum_{r\in E} e^*_r(x) \Bigr| 
           &\leq \frc{m} \biggl( \|x_i\|_{X_n} + \Bigl| 
                 \sum_{r\in E} e^*_r \sum_{j=i+1}^m x_j \Bigr| \biggr) 
            \leq \frc{m} \biggl( 1 + \sum_{l=1}^k \Bigl| 
               \sum_{r\in E_l} e^*_r \sum_{j=i+1}^m x_j \Bigr| \biggr) \\ 
           &\leq \frc{m} \biggl( 1 + \sum_{j=i+1}^m  \sum_{l=1}^k 
                  \Bigl| \sum_{r\in E_l} e^*_r(x_j) \Bigr| \biggr) 
            \leq \frc{m} \Bigl( 1 + \sum_{j=i+1}^m k\|x_j\|_{X_{n-1}} 
                                                               \Bigr) \\
           &\leq \frc{m} \Bigl( 1 + \sum_{j=i+1}^m 
                                     \frac{k_i}{k_{j-1}}2^{-j+1} \Bigr) 
            \leq \frc{m} \Bigl( 1 + \sum_{j=i+1}^m 2^{-j+1} \Bigr) 
            \leq \frac{1}{m}(1+1)=\frac{2}{m} \ ,
  \end{align*}
  and hence $\|x\|_{X_n}\leq2/m<\e$ as required.  This completes the
  proof.
\end{proof}

Lemma~\ref{lem:Xn-lli} now follows and the proof of
Theorem~\ref{thm:Xa-Ca-indices} is finished.

\end{document}